\documentclass[10pt]{article}
\usepackage{mathpazo}
\usepackage[a4paper,width=150mm,top=25mm,bottom=25mm]{geometry}
\usepackage{amsmath, amssymb, amsfonts, amsthm}
\usepackage{mathtools}
\usepackage{appendix}
\usepackage{graphicx}
\usepackage{comment}
\usepackage{hyperref}
\usepackage{titlesec}
\usepackage{authblk}

\hypersetup{
    colorlinks=true,
    linkcolor=magenta,
    citecolor=blue, 
    filecolor=magenta,      
    urlcolor=cyan,
    pdfpagemode=FullScreen,
}

\newcommand{\X}{\mathbb{X}}
\newcommand{\Z}{\mathbb{Z}}
\newcommand{\R}{\mathbb{R}}
\newcommand{\K}{\mathbb{K}}
\newcommand{\C}{\mathbb{C}}

\newcommand{\pa}{\mathrm{pa}} 
\newcommand{\anc}{\mathrm{anc}}
\newcommand{\des}{\mathrm{des}}

\newcommand{\paths}{\mathbf{P}}
\newcommand{\treks}{\mathbf{T}}
\newcommand{\G}{\mathcal{G}}
\newcommand{\D}{\mathcal{D}}
\newcommand{\E}{\mathbb{E}}
\newcommand{\F}{\mathcal{F}}

\newcommand{\fcum}{\mathcal{S}}

\newcommand{\internal}{\mathrm{I}}

\newcommand{\sh}{\mathrm{sh}}
\newcommand{\str}{\mathrm{str}}
\newcommand{\defl}{\mathrm{defl}}
\newcommand{\infl}{\mathrm{infl}}
\newcommand{\trektop}{\mathrm{top}}

\newcommand{\bmat}{\mathsf{BMat}}
\newcommand{\ebdg}{\mathsf{E}}

\newcommand{\acs}{\mathsf{ACS}}
\newcommand{\stable}{\mathsf{stb}}
\newcommand{\smry}{\mathsf{smry}}

\theoremstyle{plain}
\newtheorem{theorem}{Theorem}
\newtheorem{proposition}[theorem]{Proposition}
\newtheorem{lemma}[theorem]{Lemma}
\newtheorem{corollary}[theorem]{Corollary}

\theoremstyle{definition}
\newtheorem{definition}{Definition}
\newtheorem{example}{Example}
\newtheorem{remark}{Remark}

\newtheorem{question}{Question}

\newtheorem{convention}{Convention}

\titleformat{\section}[block]{\centering\large\bfseries}{\thesection.}{0.5em}{}[]
\titleformat{\subsection}[block]{\centering\normalsize\bfseries}{\thesubsection.}{0.3em}{}[]

\title{\scshape \LARGE $P$-stationary structural equation modeling: \\
\large Causal effect identification \\ under temporally incomplete observations}
\author{\normalsize Nicolas-Domenic Reiter \thanks{\texttt{nicolas-domenic.reiter@tum.de}}}
\affil{\normalsize Department of Mathematics, TUM School of Computation, Information and Technology, Technical
University of Munich, Munich, Germany}
\date{}

\begin{document}

\maketitle
\begin{abstract}
    Causal modeling of discrete-time processes often involves the assumption that the modeled processes are stationary. A general class of linear causal models for discrete-time stationary processes are structural equation processes (SEP), which generalize the well known structural equation model for i.i.d. observations. In this paper, we extend the SEP model from stationary to $P$-stationary processes, which are non-stationary processes whose $P$-periodically sampled subprocesses are jointly stationary. We will use this model to approach the problem of causal effect identification under temporally incomplete observations from a graphical modeling perspective. To have a notion of generic identifiability, we restrict our analysis to $P$-stationary structural vector autoregressive (SVAR) processes, which are SVAR processes whose underlying coefficients can vary $P$-periodically over time. Our findings allow to detect, in any given temporally incomplete observational setting, if a set of $P$-periodic instrumental subprocesses qualifies to identify a causal effect between $P$-periodic subprocesses, and if this is the case, that effect admits an explicit expression as a rational function evaluated at the observational spectral density. \\

\end{abstract}

\section{Introduction}
Causal inference \cite{pearl2009causality, spirtes2000causation} is about modeling cause and effect relations in the presence of randomness. The fundamental concept in causal inference is that of a causal model. A causal model has an associated directed graph whose vertices represent the entities between which causal relations are to be modeled, and whose directed links specify the immediate cause and effect relations. In addition to the graph, a causal model contains a set of functions compatible with the directed links, and a random noise object for each vertex. If the model is well defined, it gives rise to the observational distribution, which is a probability distribution over a collection of random variables indexed by the vertices of the graph. The question of what information about the causal model can be recovered from the observational distribution drives a major stream of research in causal inference, including the present paper. 

Time is intrinsic to causal relations in that an effect cannot precede its cause, i.e., cause and effect are ordered by time. When defining a causal model which obeys the ordering of time, a discrete set of time points is often assumed. Structural equation processes \cite{reiter2023causalI, reiter2026causalII} provide a broad class of linear causal models for discrete-time processes. The qualitative causal structure of an SEP is encoded by a finite directed graph. Every vertex on that graph has an independent stationary noise process, and every directed link has an associated infinite sequence, which quantifies how the past of the causing process linearly influences the present state of the caused process. The structural equation process defines a linear equation system, which is uniquely solved by the (stationary) observational process of the SEP. 

In an SEP, the influence of one process on another is time invariant, i.e., the sequence describing how the latter process at a given time point is linearly influenced by the past of the former does not depend on the time point. In particular, stationary SEPs are not capable of modeling time-varying interactions. In the current paper, we extend the class of structural equation processes to a subclass of non-stationary discrete-time processes, which we call $P$-stationary, and which allow to model periodically varying interactions between processes. For the purpose of defining $P$-stationarity, we declare a subsampling to be a subset $S \subset \Z$, and we say that it has periodicity $P$ if $t \in S$ implies that also $t \pm P \in S$. Note that there are exactly $P$-many subsamplings with periodicity $P$. In this work, a process $\X$ is called $P$-stationary if the subprocesses $\X _S =(\X(t) :{t \in S})$, where $S$ indexes the $P$-periodic subsamplings, are jointly stationary. To give an example, imagine we measure the temperature at some fixed location twice a year, once during summer and once during winter. The resulting temperature time series will turn out non-stationary if at that location winter and summer are significantly different from each other. Nonetheless, the time series consisting only of the winter and summer measurements, respectively, might still be stationary, as we would expect summer and winter, respectively, to be similar from year to year. If, in addition, the dependence between winter and summer does not change over time, the entire time series, which includes both winter and summer values, will be 2-stationary.

In this paper, we define a class of linear causal models for $P$-stationary processes whose components are indexed by some finite index set $V$. To stay in line with previous terminology \cite{reiter2023causalI, reiter2025sep, reiter2026causalII}, we term them $P$-stationary SEPs. A $P$-stationary SEP consists of a graph $G= (V\times \{0, \dots , P-1\}, D)$, a link filter for every link on $G$, and a noise process for every vertex on $G$. A vertex on $G$ represents a $P$-periodic subsampling of some component process $v$. Consequently, a link between two vertices expresses a linear influence of one $P$-periodic subprocess onto another $P$-periodic subprocess. 

Our main motivation to introduce $P$-stationary SEPs is to provide a graphical model for handling temporally incomplete observations. For the biyearly temperature measurements, the subprocess consisting only of the summer or winter values is a simple example of temporally incomplete observation. However, our notion of temporal incompleteness allows for much more general observational scenarios. 

The underlying observation to our approach is that any temporally incomplete observation of an $P$-stationary SEP can be represented, for an appropriate integer $R$, as a finite collection of $RP$-periodic subsamplings. Then we will see that every $P$-stationary SEP has an equivalent $RP$-stationary SEP, which is equipped with a latent structure rendering every missing subsample as latent. In this setting, we analyze the problem of identifying causal effects, i.e., link filters, in $P$-stationary SEPs from temporally incomplete observations. Following previous work \cite{reiter2025sep, reiter2026causalII}, we will study the problem in the frequency domain, where the algebraic relations are more tractable. The identification of a given effect in the $P$-stationary SEP using a set of instrumental processes requires certain entries in the spectral density. If the observations are temporally incomplete, one faces the problem that the spectral density can generally not be computed, so identification via instrumental processes seems infeasible. However, by providing a linear recovery formula (Proposition \ref{prop: deflation in the frequency domain}), we can characterize (Definition \ref{def: linear identifiability spectrum}) the temporally incomplete observational settings for which it is possible to recover a given entry in the spectral density (using that linear formula). Accordingly, in a setting of temporally incomplete observations, we detect a given effect to be identifiable if there is a suitable set of instrumental processes while the temporal incompleteness of the observations is such that all the entries in the spectral density that are required for identification via that set of instrumental processes can be recovered (Theorem \ref{thm: instrumental process identification}). For technical reasons, we restrict our analysis to $P$-stationary structural vector autoregressive processes, i.e., SVAR processes whose coefficients are allowed to vary periodically over time with periodicity $P$.

To exemplify the problem and our way of approaching it, consider the simple setup of three stationary SVAR processes, whose causal relations are given by $u \to v \to w$. Suppose we are interested in the effect associated with the link $v \to w$. If we had temporally complete observations for $v$ and $w$, we could identify this effect from the spectral density via linear regression in the frequency domain. However, this approach cannot be pursued if the observations for $v$ and $w$ are temporally incomplete, since in this case the spectral density cannot be computed. Assuming temporally complete observations for the instrumental process $u$, it will follow from the theory we develop in this paper that the spectral dependence between $u$ and $v$, and between $u$ and $w$ can be recovered, so that the effect $v \to w$ can be identified by the usual instrumental process formula for the frequency domain. A more interesting case, discussed in Example \ref{example: iv identification 2} and pictured in Figure \ref{fig: example 2 iv identification}, demonstrates that our findings allow to identify effects even if none of the involved processes is fully observable.

The rest of the paper is organized as follows. In Section \ref{section: preliminaries}, we recall some basic terminology for linear causal models of discrete-time stationary processes. In Section \ref{section: P-stationary processes}, we define $P$-stationary SEPs and SVAR processes and show how an $P$-stationary SVAR process can equivalently by described as an $P$-stationary SEP. Furthermore, we define the embedding operator that identifies for every $P$-stationary SEP an equivalent $RP$-stationary SEP. In Section \ref{section: frequency domain}, we present the algebraic structure of $P$-stationary SEPs in the frequency domain. The main result of this section will be used to determine if a given temporally incomplete observational pattern allows to linearly recover a specified entry in the spectral density. In Section \ref{section: identification}, we formalize the problem of identifying causal effects in $P$-stationary SEPs from temporally incomplete data in the frequency domain, i.e., we will make precise what we mean by a temporally incomplete observations and we will provide a notion of rational identifiability for $P$-stationary SVAR processes (under a given temporally incomplete observational pattern). Using those terms, we will formulate the identifiability result. 

\subsection{Related work} 
We approach causal inference in the framework of \emph{Structural Causal Models} presented and developed in, e.g., \cite{pearl2009causality, hoyer2008nonlinear, buhlmann2014cam, spirtes2000causation}. Since we focus on linear models, much of the work and theory developed here is inspired and influenced by the theory and terminology used in the literature on linear \emph{Structural Equation Modeling}, see, e.g., \cite{bollen1989structural, shimizu2006linear, wright1921correlation, wright1934method}, and especially those contributions that take an algebraic perspective on the matter, see, e.g., \cite{drton2007algebraic, drton2009lectures, foygel2012halftrek, drton2011global, sullivant2010trek}. 

Our paper is particularly related to those contributions which incorporate time as a discrete index into the modeling of causal relations. Early contributions include those related to \emph{Granger Causality}, see \cite{geweke1984inference, geweke1984measures, granger1963economic, granger1969investigating, granger1980testing}, which underwent a reformulation in terms of concepts and terminology from graphical modeling, see \cite{dahlhaus2000graphical, dahlhaus2003causality, eichler2005graphical, eichler2007causal, eichler2007granger, eichler2012graphical, eichler2013causal}. This line of work has then been further developed, for instance, in \cite{runge2019detecting, runge2019inferring, runge2020discovering, gerhardus2020highrecall, hyvarinen2010estimation}. 

The theory presented in the present paper builds on \cite{reiter2023causalI, reiter2024asymptoticuncertaintyestimationfrequency, reiter2025sep}, which develop the notion of (stationary) structural equation processes to quantify process-level causal relations (in the frequency domain). An algebraic approach to structure and effect identification from observational data, in the spirit of the contributions mentioned above, was developed in \cite{reiter2026causalII}.

A specific from of temporally incomplete observations are undersampled time series, i.e., when observations are gathered less frequently than causal effects happen. An introduction to the problems related to inferring causal relations from undersampled time series can be found in \cite{danks2013learning}. Specific approaches to inferring the causal structure from an undersampled time series model have been developed in, e.g., \cite{liu2023causal, hyttinen2016causal}. The special case where the underlying time series model is linear is investigated, for example, in \cite{gong2015discovering, tank2019mixedfrequency}. Specifically, they consider SVAR processes of order one (causal effects happen either instantaneously or with a time lag of one time unit) generated by a non-Gaussian white noise process. In \cite{tank2019mixedfrequency}, the observations are allowed to be gathered at two different sampling rates, which is also known as the mixed-frequency case. 
Another approach to the identification of VAR models (without instantaneous effects) is presented in \cite{anderson2012identifiability, anderson2016multivariate}, where the authors consider stable VAR processes of some arbitrary finite order $p$, and which are assumed to be sampled at two different frequencies. The fast frequency is assumed to be the actual time scale of the process, i.e., the frequency with which causal effects happen. Under the assumption that at least one process is observed at the actual time scale, they show that for almost every possible VAR parameter configuration, that parameter can be uniquely recovered from the second moments of the observed process (which is generally temporally incomplete).

In this work, we consider linear time series models whose coefficients are allowed to vary periodically over time. Such models generally yield non-stationary processes. What we call an $P$-stationary SVAR process in this work is a special instance of a time varying VAR process, see, e.g., \cite{subbarao1970fitting, dahlhaus1997fitting}. A recent contribution to non-stationary undirected graphical modeling can be found in \cite{basu2023graphical}. Our approach to weaken the stationarity assumption in a discrete-time setting is similar to the one taken in \cite{gao2023causal}, which develops a causal discovery algorithm for so-called semi-stationary Markov chains with an underlying periodically varying causal structure.  
\subsection{Basic notation}
Let $V$ be some finite set and $\K$ a field, then $\K^V$ denotes the $\K$-vector space of dimension $|V|$ and whose coordinates are indexed by $V$. Let $U$ be another finite set, then $\K^{V \times U}$ denotes the space of matrices whose rows are indexed by $V$ and whose columns are indexed by $U$. In case $V=U$, we write $I_V$ for the identity matrix. For a matrix $A=[A_{v,u}]_{v \in V, u \in U} \in \K^{V \times U}$, and subsets $V' \subset  V$ and $U' \subset U$, we define the submatrix of $A$ corresponding to the entries $V' \times U'$ as $[A]_{V', U'}\coloneqq [A_{v,u}]_{v \in V', u \in U'} \in \K^{V' \times U'}$. The transpose of $A$ is denoted $A^\top \coloneqq[A_{v,u}]_{u \in U, v \in V} \in \K^{U \times V}$. Since we will often encounter block-matrices, it will be useful to have a dedicated notation for them. Let $U,V,I,J$ be finite sets, then define $\bmat(\K, (U,V), (I,J))$ to be the set of matrices $A$ which have the following two level block structure: 
\begin{align*}
    A & = \begin{bmatrix}
        A_{u,v}
    \end{bmatrix}_{u \in U,v \in V}, & A_{u,v} &= \begin{bmatrix}
        A_{u,v}^{(i,j)}
    \end{bmatrix}_{i\in I,j \in J } \in \K^{I \times J}.
\end{align*}
If $U=V$ or $I= J$, we will abbreviate the corresponding set of block-matrices by $\bmat(\K, V, (I,J))$ and $\bmat(\K, (U,V), I)$ respectively, and if both $U =V$ and $I= J$, then we abbreviate the corresponding set of block-matrices by $\bmat(\K, U, I)$.  

For two integers $p,q \in \Z$, define $[p,q] \coloneqq \{ i \mid p \leq i \leq q\}$, and  for $n \geq 1$, we write $[n]$ to abbreviate $[1,n]$. For a positive integer $P \geq 1$, we define the set of $P$-multiples to be $P \Z  \coloneqq \{kP \mid k \in \Z\}$, and for $m \in [0, P-1]$, we define $ P\Z -m\coloneqq \{ kP-m \mid k \in \Z \}$. 


For $\K$ either $\R$ or $\C$, the space $\ell^p(\Z, \K)$ consists of all sequences $\lambda: \Z \to \K$ such that $\lVert \lambda \rVert_p^p \coloneqq \sum_{k \in \Z } |\lambda(k)|^p < \infty$. 
We use $\ell^p(\Z, \K^{V \times U})$ to denote the space of sequences $\Lambda: \Z \to \K^{V \times U}$ such that $[\Lambda]_{v,u} \coloneqq ( k\in \Z \mapsto [\Lambda(k)]_{v,u}) \in \ell^p(\Z, \K)$ for all $(v,u) \in V \times U$. In this work, we will call any element in $\Lambda \in \ell^1(\Z, \K^{V\times U})$ such that $\Lambda(k)= 0$ if $k<0$, a \emph{filter}, a term often used in the time series literature. For two sequences $\lambda, \mu\in \ell^p(\Z,\K)$, the \emph{convolution} resp. \emph{tilted convolution} are defined as follows:
\begin{align*}
    \lambda \ast \mu &\coloneq \left( k \mapsto \sum_{l \in \Z } \lambda(l)\mu(k-l) \right), & \lambda \ast' \mu &\coloneqq \left( k \mapsto \sum_{l \in \Z}\lambda(k+l)\mu(l) \right).
\end{align*}
The \emph{transpose} of a sequence $\Lambda$ is the element-wise transpose, $\Lambda^\top \coloneqq (k \mapsto \Lambda(k)^\top)$. If $V = U$, then $\mathbb{I}_V$ denotes the \emph{identity filter}, i.e., $\mathbb{I}_V(0) = I_V$, and $\mathbb{I}_{V}(k) = 0$ for all $k \neq 0$. Furthermore, for a fixed set $S$ and a subset $S' \subset S$, we write $\mathbf{1}_{S'}:S \to \{0,1\}$ for the \emph{indicator function}, i.e., $\mathbf{1}_{S'}(x) = 1$ if and only if $x \in S'$.  

\section{Preliminaries}\label{section: preliminaries}
In this section, we review the basic terminology for directed graphs, discrete-time stationary processes, and linear causal models of processes, i.e., structural equation processes \cite{reiter2023causalI, reiter2025sep, reiter2026causalII}. 
\subsection{Directed graphs}
A directed graph $G = (V,D)$ consists of a countable set (finite or infinite) of vertices $V$ and a set of directed links $D \subset V \times V \setminus \{ (v,v) \mid v \in V\}$. To express that $(v,w) \in D$ is a directed link, we write $v \to w$ for short, and $v \not \to w$ otherwise. Note that we exclude self-links, i.e., links of the form $ v \to v$. The \emph{parent set} of $v \in V$ is defined to be $\pa(v) \coloneqq \{ u \in V | u \to v \}$. A \emph{path} on $G$ is a finite sequence $\pi=(v_i \to w_i)_{i \in [n]}$ of links such that either $w_i = v_{i+1}$ or $w_i = w_{i+1}$ or $v_i = v_{i+1}$ for all $1 \leq i < n$. A path $\pi$ is called \emph{directed} if $w_i = v_{i+1}$ for all $1\leq i < n$. By convention, we also consider empty paths, i.e., for every vertex $v \in V$, we denote by $\epsilon_v$ the path from $v$ to itself that does not contain any link. For two vertices $v,w \in V$, we will denote the set of all directed paths from $v$ to $w$ by $\paths_G(v,w)$. A graph $G$ is called \emph{acyclic} if $\paths_G(v,v) = \{ \epsilon_v \}$ for every $v \in V$. The set of \emph{ancestors} and \emph{descendants}, respectively, of a vertex $v \in V$ are defined as follows:
\begin{align*}
    \anc_G(v) &\coloneqq \{ u \in V \mid \paths_G(u,v) \neq \emptyset\} & \des_G (v) & \coloneqq \{ w \in V | \paths_G(v,w) \neq \emptyset\}. 
\end{align*}
The \emph{strongly connected component} of a vertex $v \in V$ is the intersection $\anc_G(v) \cap \des_G(v)$. Note that the vertex set of any directed graph can be partitioned into its strongly connected components.

Given a subset $U \subset V$, we write $\pa_G(U) = \bigcup_{u \in U } \pa_G(u)$ to refer to the parents of the vertices in $U$, and analogously for the ancestors $\anc_G(U)$ and the descendants $\des_G(U)$. 

A \emph{trek} between two vertices $v,w \in V $ on $G$ is a tuple $\tau = (\tau_L, \trektop(\tau), \tau_R)$, where $\trektop(\tau) \in V$ and $\tau_L \in \paths_G(\trektop(\tau), v)$ and $\tau_R \in \paths_G(\trektop(\tau), w)$. The set of treks between two vertices $v,w \in V$ on $G$ is denoted by $\treks_G(v,w)$. 

For a subset $U \subseteq V$, the restriction of $G$ to $U$ is the directed graph $G_U = (U, D_U)$, where $D_U = D \cap U \times U$. Furthermore, suppose $(G_i = (V_i, D_i))_{i \in I }$ is a discrete collection of directed graphs, then $\bigcup_{i \in I}G_i = (\bigcup_{i \in I }V_i, \bigcup_{i \in I} D_i )$ is the union graph and, similarly, $\bigcap_{i \in I} G_i = (\bigcap_{i \in I} V_i, \bigcap_{i \in I} D_i)$ is the intersection graph. 

Given a finite set $V$, a directed graph $\G=(V \times \Z, \D)$ is called a \emph{(directed) time series graph} if for any link $v(s) \to w(t)$, it holds that $s \leq t$, i.e., links do not point backwards in time. A time series graph $\G$ is said to be \emph{stationary} if $v(t-k) \to w (t)$, for some $t\in \Z$, implies $v(s-k) \to w(s)$ for all $s \in \Z$. Given a stationary time series graph $\G$, we write $v \to_k w$ to abbreviate that $v(t-k) \to w(t)$ for every $t \in \Z$. A stationary time series graph $\G=(V \times \Z, \D)$ can be reduced to a finite graph $G=(V,D)$ that contains the directed link $v \to w$, where $v \neq w$, if and only if and there is some $k \geq 0$ such that $v \to_k w$ holds on $\G$. Following \cite{reiter2023causalI, reiter2025sep, reiter2026causalII}, we call this finite reduction the \emph{process graph} of $\G$. 

\subsection{Stationary discrete-time processes}
Let $V$ be a finite set and let $\X = (\X_v(t)\mid {v \in V, t \in \Z})$ be a $V$-variate discrete-time stochastic process. In this work, we use the word process to refer to a discrete-time stochastic process. For subsets $U \subseteq V$ and $T \subseteq \Z$, we define the associated subprocess $\X_U(T) \coloneqq  (\X_u(t)\mid {u \in U, t \in T })$. For a subset $\mathbf{t} \subset \Z$ and some $s \in \Z$, we define $\mathbf{t} - s \coloneqq \{ t-s \mid t \in \mathbf{t}\}$.  A process $\X$ is called \emph{stationary} if for every finite $\mathbf{t} \subset \Z$ and every $s\in \Z$, the random vectors $\X(\mathbf{t})$ and $\X(\mathbf{t} - s )$ follow the same distribution. Otherwise, the process $\X$ is called \emph{non-stationary}. For pairwise disjoint subsets $U_1, \dots, U_n \subset V $, the associated subprocesses $\X_{U_1}, \dots, \X_{U_n}$ are called mutually independent if for any finite $\mathbf{t}\subset \Z$ the random vectors $\X_{U_1}(\mathbf{t}), \dots ,\X_{U_n}(\mathbf{t})$ are mutually independent. The process $\X$ is called a \emph{white noise} process if for any two finite subsets $\mathbf{t}, \mathbf{s} \subset \Z$ such that $\mathbf{t} \cap \mathbf{s} = \emptyset$, the random vectors $\X(\mathbf{t})$ and $\X(\mathbf{s})$ are independent. The mean of a stationary process $\X$ is $\mu \coloneqq \E[\X(t)]$ and the \emph{autocovariance sequence (ACS)} is given by $\acs(\X) \coloneqq (k \mapsto  \E [\X (t) \X(t-k)^\top])$. 

\subsection{Structural Equation Processes}
A \emph{structural equation process (SEP)} $\mathbb{S}=(G,\Lambda, \X^\internal)$ consists of: a finite directed graph $G=(V,D)$; a matrix-valued filter $\Lambda \in \ell^1(\Z, \R^{V \times V})$ such that $v \not \to w$ implies that $\Lambda_{v,w}$ is the zero-filter; and a $V$-variate mutually independent process $\X^\internal$. A solution to $\mathbb{S}$ is a process $\X$ such that it satisfies $\X = \Lambda^\top \ast \X + \X^\internal$.
An SEP $\mathbb{S}$ has a (unique) solution if the filter $\Lambda$ is such that $S(\Lambda) \coloneqq \sum_{k \geq 0} \Lambda^k$ exists and is absolutely summable. In this case, the unique solution is given by the process $S(\Lambda^\top) \ast \X^\internal$. 

One way to construct an SEP is to start with an underlying \emph{structural vector autoregressive (SVAR)} process, which consists of a stationary time series graph $\G= (V \times \Z, \D)$; a filter $\Phi: \Z \to \R^{V \times V}$ such that $ v \not\to_k w$ implies $\Phi_{v,w}(k) = 0$; and a $V$-variate mutually independent white noise process $\eta $. A stationary process $\X$ is said to be an SVAR process compatible with $(\G, \Phi, \eta)$ if it satisfies $\X = \Phi^\top \ast \X  + \eta$.
Such a process $\X$ exists uniquely (up to a zero-set) if the coefficients $\Phi$ are \emph{stable}, see \cite{luetkepohl2005timeseries, reiter2026causalII}. Additionally, every stable (in the sense of \cite{reiter2026causalII}) SVAR process can be represented as an SEP $\mathbb{S}=(G, \Lambda, \X^\internal)$, where $G$ is the process graph of $\G $, see \cite{reiter2023causalI, reiter2025sep, reiter2026causalII}, meaning that the SVAR- and SEP solution coincide. 

\section{A structural equation model for $P$-stationary processes}\label{section: P-stationary processes}
We begin this section by introducing an embedding operator for sequences and processes. This operator will be central to the development of the main results of this paper, since it will allow us to handle temporally incomplete observations. Using the embedding operator, we introduce $P$-stationary processes and filters to define $P$-stationary structural equation processes as a general class of linear causal models for $P$-stationary processes. In the final part of this section, we will see that this class encompasses structural vector autoregressive processes whose coefficients can vary $P$-periodically over time. This observation will be important, since in developing the main results (stated for such $P$-stationary SVAR processes) we will use the structure of an $P$-stationary SEP as well as the additional structure of an $P$-stationary SVAR processes.
\subsection{An embedding operator}
To prepare the construction of the embedding operator, we introduce the following three operators on the space of filters. 
Let $\lambda: \Z \to \K^{V\times U}$ be a sequence. For $j \in \Z$, the \emph{$j$-shifted} sequence is defined as $\sh^j(\lambda) \coloneqq \left(k  \mapsto \lambda(k-j) \right)$. For a positive integer $N$, we have the \emph{$N$-deflated} sequence $\defl_N(\lambda) \coloneqq (k \mapsto \lambda(kN))$ and the \emph{$N$-inflated} sequence $\infl_N(\lambda) \coloneqq (k \mapsto \mathbf{1}_{N\Z}(k) \lambda(\mathbf{1}_{N\Z}(k)kN^{-1}))$. Based on these operators, we define for a positive integer $N$ two linear embedding operators
\begin{align*}
    \mathsf{e}_N, \mathsf{e}_N' : \ell^2( \C^{U \times V }) \to \ell^2 ( \bmat(\C, [0,N-1], (U,V))). 
\end{align*}
Those operators map a sequence $\lambda \in \ell^2 ( \C^{U \times V})$ to the following sequences:
\begin{align*}
     \mathsf{e}_N(\lambda ) &\coloneqq \begin{bmatrix}
        \defl_N \circ \sh^{i-j} (\lambda) 
     \end{bmatrix}_{i,j \in [0,N-1]}, & \mathsf{e}_N' (\lambda) &\coloneqq \begin{bmatrix}
          \defl_N \circ \sh^{j-i} (\lambda )
     \end{bmatrix}_{i,j \in [0,N-1]}. 
\end{align*}
Analogously, we define for a $V$-variate process $\X$, the embedded process $\mathsf{e}_N(\X) \coloneqq [\defl_N \circ \sh^{i} (\X)]_{i \in [0,N-1]}$, which is a process that consists of $N|V|$ component processes.  
In view of the theory we develop in this work, the following two lemmas will be important.
\begin{lemma}\label{lemma: linear recovery of the embedded filter}
Let $\lambda \in \ell^2( \C^{U \times V })$, then for every $i\in [0,N-1]$ it holds that
\begin{align}\label{shifted recovery property}
    \sh^{i}(\lambda) &= \sum_{j=0}^{N-1} \sh^{j} \circ \infl_N\left( \defl_N \circ \sh^{i-j} (\lambda)\right).
\end{align}
\end{lemma}
\begin{lemma}\label{lemma: (tilted) convolution under the embedding}
    Let $\lambda \in \ell^2(\C^{W \times V} )$ and $\mu \in \ell^2(\C^{V \times U})$, then it holds that 
    \begin{align*}
        \mathsf{e}_N(\lambda \ast \mu ) &= \mathsf{e}_N(\lambda) \ast \mathsf{e}_N(\mu ), & \mathsf{e}_N (\lambda \ast ' \mu ) &= \mathsf{e}_N (\lambda) \ast' \mathsf{e}_N'(\mu ). 
    \end{align*}
    Furthermore, for a $V$-variate process $\X$, it holds that $\mathsf{e}_N (\lambda \ast\X) = \mathsf{e}_N(\lambda) \ast \mathsf{e}_N(\X) $.
\end{lemma}

\subsection{$P$-stationary processes}
\begin{definition}[$P$-stationary process]\label{def: P-stationary process}
    Let $P$ be a positive integer. A process $\X = (\X_v \mid {v \in V})$ is called \emph{$P$-stationary} if the stacked process $(\mathsf{e}_P(\X_v)\mid v \in V)$ consisting of $|V|P$ many components, is stationary. For two distinct components $u ,v \in V$, the associated component processes are said to independent if the stationary processes $\mathsf{e}_P(\X_v)$ and $\mathsf{e}_P (\X_u)$ are independent. 
\end{definition}
\begin{remark}\label{remark: P-stationarity => kP-stationarity}
    Given a positive integer $R$, every $P$-stationary process is also $RP$-stationary. 
\end{remark}
\begin{convention}\label{p-stationary filtering as block-wise filtering}
    In this work, we will implicitly consider any $V$-variate $P$-stationary process $\X$ in its stationary from. That means, when introducing $\X$ as $P$-stationary, we will formally refer by $\X$ to the stationary process $(\mathsf{e}_P(\X_v) \mid v \in V)$. 
\end{convention}
\begin{definition}[$P$-stationary ACS]
    Let $\X$ be a $V$-variate $P$-stationary process.
    The \emph{$P$-stationary autocovariance sequence (ACS)} of $\X$ is the ACS of the stationary process $(\mathsf{e}_P(\X_v) \mid v \in V)$. In particular, the $P$-stationary ACS of $\X$ is block-structured such that $\acs(\X): \Z \to \bmat(\R, V, [0,P-1])$. 
\end{definition} 
\begin{definition}[Filtering $P$-stationary processes]\label{def: P-stationary filtering}
    Let $V$ and $W$ be finite sets, a block-structured filter of the form $\Lambda: \Z \to \bmat(\R, (V,W), [0,P-1])$ is said to be a $P$-stationary filter. Given a $V$-variate $P$-stationary process $\X$, the process $\Lambda^\top \ast \X$ is said to be the \emph{$\Lambda$-filtering} of $\X$.
\end{definition} 
\begin{lemma}\label{lemma: acs p-stationary filtering}
    Let $\X$ be a $V$-variate $P$-stationary process, and $\Lambda: \Z \to \bmat(\R, (V,W), [0,P-1])$ a $P$-stationary filter, then it holds that $\acs(\Lambda^\top \ast \X) = \Lambda^\top \ast \acs(\X) \ast' \Lambda$.
\end{lemma}
\begin{question}\label{question: consistent embedding of P-stationary filters}
    As of Remark \ref{remark: P-stationarity => kP-stationarity}, for a $P$-stationary process $\X$ and $P$-stationary filter $\Lambda$, the processes $\X$ and $\mathbb{Y}\coloneqq\Lambda \ast \X$ are also $RP$-stationary. How to describe the linear relation between $\X$ and $\mathbb{Y}$ in terms of an $RP$-stationary filter?
\end{question}
To answer this question, we make use of the embedding operator constructed in the previous subsection. Let $\Lambda$ be a $P$-stationary filter, then we define the block-wise embeddings as follows: 
\begin{align*}
    \ebdg_R (\Lambda) &\coloneqq \begin{bmatrix}
        \mathsf{e}_R (\Lambda_{v,w})
    \end{bmatrix}_{v,w \in V }, & \ebdg_R'(\Lambda) &\coloneqq  \begin{bmatrix}
        \mathsf{e}_R' (\Lambda_{v,w})
    \end{bmatrix}_{v,w \in V}.
\end{align*}
Note that $\mathsf{E}_{R}(\Lambda)$ has a three-level block structure, where the first level is indexed by $V$, the second by $[0,R-1]$ and the third by $[0,P-1]$.
Given a $V$-variate $P$-stationary process $\X$ in its stationary form (see Convention \ref{p-stationary filtering as block-wise filtering}), the embedding applied to $\X$ takes the following stacked form:
    \begin{align*}
        \ebdg_{R} (\X) &= [\mathsf{e}_{R}\circ \mathsf{e}_P(\X_v)]_{v \in V}.
    \end{align*}
In particular, $\mathsf{E}_{R}(\X)$ is the stationary form, see Convention \ref{p-stationary filtering as block-wise filtering}, of $\X$ read as an $RP$-stationary process. The following states that the block-wise embedding is compatible with convolution, filtering of processes and the computation of the ACS. 
\begin{proposition}\label{prop: embedding properties}
    Let $\Lambda : \Z \to \bmat(\R, (V,W)), [0,P-1])$ and $\Gamma: \Z \to \bmat(\R, (U,V), [0,P-1])$ be two $P$-stationary filters, then it holds that
    \begin{align*}
        \mathsf{E}_{R}(\Gamma \ast \Lambda) &= \mathsf{E}_{R}(\Gamma) \ast \mathsf{E}_{R}(\Lambda), & \mathsf{E}_{R}(\Gamma \ast ' \Lambda)&= \mathsf{E}_{R}(\Gamma) \ast' \mathsf{E}'_{R}(\Lambda).  
    \end{align*}
    Furthermore, for a $V$-variate $P$-stationary process $\X$, it holds that $\mathsf{E}_{R}(\Lambda^\top \ast \X) = \mathsf{E}_{R}(\Lambda^\top) \ast \mathsf{E}_{R}(\X)$, and $\acs(\mathsf{E}_{R}(\X)) = \mathsf{E}_{R}(\acs(\X))$. 
\end{proposition}
    
\subsection{$P$-stationary structural equation processes}
For a finite set $X$, we call a filter $\Lambda : \Z \to \R^{X \times X}$ compatible with a finite directed graph $G= (X ,D)$ if $v \not \to w $ implies that $\Lambda_{v,w}$ is the zero filter. 
\begin{figure}
    \centering
    \includegraphics[width=0.8\linewidth]{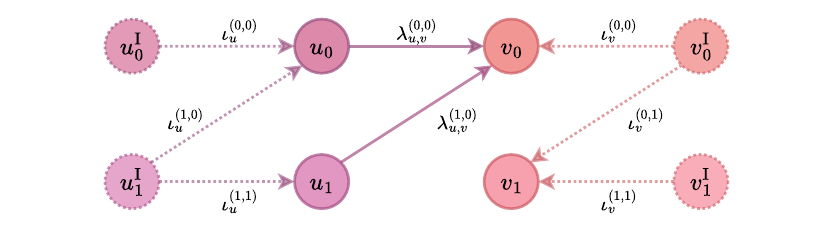}
    \caption{For $V= \{u,v\}$, the graph specifies a $2$-stationary SEP together with an IDFS. A vertex of the form $a_k$ and $a_k^\internal$, respectively, represents the subprocess $\defl_2 \circ \sh^k(\X_a)$ and $\defl_2 \circ \sh^k (\eta_a)$, respectively.}
    \label{fig: 2-stationary sep general case}
\end{figure}
\begin{definition}[$P$-stationary SEP]
    Let $P$ be a positive integer. A \emph{$P$-stationary structural equation process (SEP)} over some finite index set $V$ is given by a triple $\mathbb{S} = (G, \Lambda, \X^{\internal})$ consisting of: a finite directed graph $G = (V \times [0,P-1], D)$; a $P$-stationary filter $\Lambda: \Z \to \bmat(\R, V, [0,P-1])$ compatible with $G$; and a zero-mean $V$-variate $P$-stationary process $\X^{\internal}$. A $P$-stationary process $\X$ is said to be a solution for the SEP $\mathbb{S}$ if it satisfies the structural equation $\X = \Lambda^\top \ast \X + \X^\internal$. 
    The $P$-stationary SEP $\mathbb{S}$ is called \emph{stable} if the filter $S(\Lambda) = \sum_{k \geq 0} \Lambda^k$ is absolutely summable.
\end{definition}  
For a stable $P$-stationary SEP $\mathbb{S}=(G, \Lambda, \X^\internal)$, the process $\X = S(\Lambda^\top) \ast \X^\internal$ is the solution to the SEP $\mathbb{S}$. In this case, Lemma \ref{lemma: acs p-stationary filtering} implies that the $P$-stationary ACS of the solution $\X$ is given by
\begin{align}\label{equation: acs p-stationary sep}
    \acs(\X) &=S(\Lambda)^\top \ast \acs(\X^\internal) \ast' S(\Lambda). 
\end{align}
So far, the internal dynamics have been taken to be some general $P$-stationary process. The simplest instance of an $P$-stationary processes is an $P$-stationary white noise process.   
\begin{definition}[$P$-stationary white noise]
    A $P$-stationary univariate process $\eta$ is called a white noise process if the stationary process $\mathsf{e}_P (\eta)$ is a white noise process such that all $P$ components are mutually independent. 
\end{definition}
In this work, the internal dynamics will always arise as filtered white noise processes. In the following, we define their general structure. For that purpose, we will use $V^\internal = \{v^\internal \mid v \in V\} $ to denote a (disjoint) a copy of $V$, i.e., $V^\internal \cap V = \emptyset$.
\begin{definition}[IDFS]
     A $V$-variate and $P$-stationary \emph{internal dynamic filtering structure (IDFS)} is given by a $P$-stationary SEP $\mathbb{S}^\internal =(G^\internal, \iota, \eta)$, where $G^\internal = (V^\internal \cup V \times [0,P-1], D^\internal)$ is such that every link is of the form $(v^\internal, m) \to (v,n)$ for some $v \in V$ and some $m,n \in [0,P-1]$; a $P$-stationary filter $\iota: \Z \to \bmat(\R, (V^\internal, V), [0,P-1])$ compatible with $G^\internal$; and a $V^\internal$-variate $P$-stationary mutually independent zero-mean white noise process $\eta$.  
     The solution process, $\X^\internal = \iota^\top \ast  \eta$,  will be referred to as the \emph{internal dynamics}. 
\end{definition}
Note that the assumption on $G^\internal$ implies that the components of $\X^\internal$ are mutually independent. In Figure \ref{fig: 2-stationary sep general case}, we display the graphical data of an $2$-stationary SEP in conjunction with that of an IDFS.    
\subsection{$P$-stationary SVAR processes}\label{section: P-stationary SVAR processes}
We will now introduce $P$-stationary time series graphs and $P$-stationary SVAR processes. Then we will construct for a given $P$-stationary SVAR process its equivalent $P$-stationary SEP representation.  
\begin{definition}[$P$-stationary time series graph]
    A time series graph $\G = (V \times \Z, \D)$ is called $P$-stationary if for any directed link $v(t-k) \to w(t)$ on $\G$, it holds that for every $l \in \Z$ also $v(t-k + lP)\to w(t+lP)$ is a directed link on $\G$. 
\end{definition}
Given a $P$-stationary time series graph $\G=(V\times \Z, \D)$, we define the following sets of links: 
\begin{align*}
    \D^{(i,j)} &\coloneqq \{ v(t-k) \to w(t) \mid t \in P\Z-j, \quad t-k \in P\Z-i \}, & i,j\in [0,P-1].
\end{align*}
It follows from $\G$ being $P$-stationary that these subsets are disjoint and cover the entire set $\D$. For a link $v(t-k) \to w(t) \in \D^{(i,j)}$, we write $v \to_l^{(i,j)}w$, where $l$ is such that $k =i-j + lP$. The order of a $P$-stationary time series graph is defined to be 
\begin{align*}
    \max_{u,v \in V, i,j\in [0,P-1]} \{ l \mid u \to_l^{(i,j)} v\}.  
\end{align*}
Throughout this work, we assume the order of any time series graph to be finite. We call a $P$-stationary filter $\Phi: \Z \to \bmat(\R, V, [0,P-1])$ compatible with a $P$-stationary time series graph $\G$ if $u \not\to_k^{(i,j)}v$ implies that $\Phi^{(i,j)}_{u,v} = 0$.  


\begin{definition}[P-stationary SVAR process]
    A $P$-stationary structural vector autoregressive (SVAR) model consists of: a $P$-stationary time series graph $\G=(V \times \Z, \D)$; a $P$-stationary filter $\Phi: \Z \to \bmat(\R, V, [0,P-1])$ compatible with $\G$; and a zero-mean $P$-stationary mutually independent white noise process $\eta$. A $P$-stationary process $\X$ solves the SVAR model if it satisfies $\X = \Phi^\top \ast \X  + \eta$.
\end{definition}
In the following, $\R^\G$ will be the vector space of filters $\Phi: \Z \to \bmat(\R, V, [0,P-1])$ compatible with the time series graph $\G$. Note that $\R^\G$ is finite dimensional, since we assume $\G$ to have finite order.
\begin{figure}
    \centering
    \includegraphics[width=\linewidth]{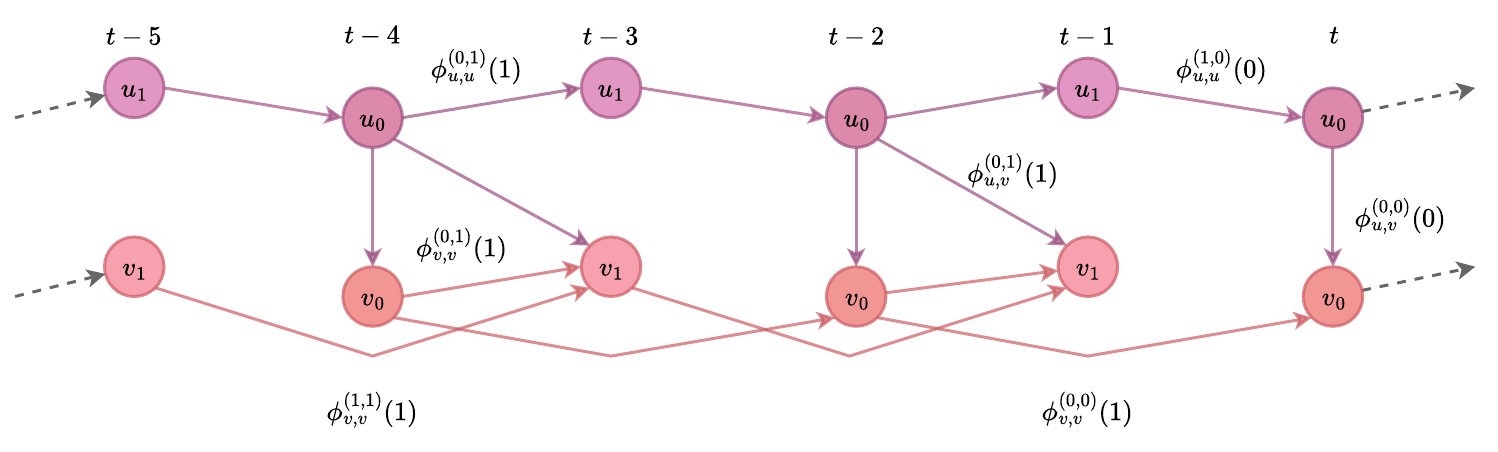}
    \caption{For $t \in 2\Z$, the graph shows an excerpt of a $2$-stationary time series graph.}
    \label{fig: 2 stationary ts graph}
\end{figure}

\begin{example}
    Consider the finite excerpt of the graph $\G = (\{u,v\}\times \Z , \D)$ depicted in Figure \ref{fig: 2 stationary ts graph} together with a filter $\Phi$ compatible with $\G$. Following the notation introduced above, the $2$-stationary graph $\G$ contains the following relations, each with an associated entry in $\Phi$. The auto-dependencies of the process $u$ are given as follows:  
    \begin{align*}
        u\to_1^{(0,1)} u: \quad  \phi_{u,u}^{(0,1)}(1),& &
        u \to_0^{(1,0)} u: \quad \phi^{(1,0)}_{u,u}(0). 
    \end{align*}
    Similarly, one can read off the relations associated with the auto-causal links of $v$ and the links arising from the lagged $2$-stationary causal links from $u$ on $v$. 
\end{example} 

In the remainder of this subsection, we will reformulate a given $P$-stationary SVAR process as an equivalent $P$-stationary SEP, meaning that the solution processes of the SVAR model and the associated SEP will coincide. For the remainder of this section, we fix an $P$-stationary SVAR model $(\G, \Phi, \eta )$.

We need to define two graphs. Namely, the graph $G$ defining the $P$-stationary SEP and the graph $G^\internal$ defining the internal dynamics. We will derive both of them from the following. 
\begin{definition}[$P$-stationary summary graph]\label{definition: summary graph}
Let $\G$ be a $P$-stationary time series graph. Its \emph{summary graph}, denoted $\smry(\G )$, is a directed graph over the vertex set $V \times [0,P-1]$ that connects two distinct vertices by $(v,i) \to (w,j)$ if and only if there is a $k\geq0$ such that $v\to_k^{(i,j)}w$ holds on $\G$.
\end{definition}

To define the structure for the internal dynamics of $v\in V$, let $\G_v$ be the restriction of $\G$ to $\{v \} \times \Z$. The summary graphs $(\smry(\G_v)\mid v \in V)$ determine $G^\internal$ in the following way:  
\begin{align*}
    \paths_{\smry(\G_v)}((v,i), (v,j)) \neq \emptyset \iff (v^\internal , i) \to (v,j).
\end{align*}
Finally, the process graph $G$ of the time series graph $\G$ is defined as follows: for $v \neq w$, there is a directed link $(v,i) \to (w,j)$ on $G$ if and only if there is some $l \in [0,P-1]$ such that: $(v,i) \to (w,l)$ is a directed link on $\smry(\G)$, and $( w^\internal, l)\to (w,j)$ is a link on $G^\internal$. 

\begin{figure}
    \centering
    \includegraphics[width=\linewidth]{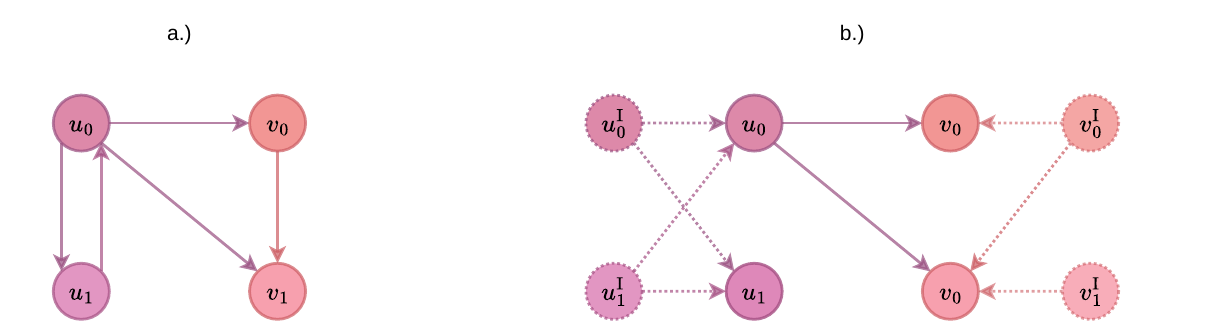}
    \caption{a.) The summary graph derived from the graph in Figure \ref{fig: 2 stationary ts graph}. b.) The graph $G \cup G^\internal$ derived from a.). In both graphs, a vertex of the form $a_k$ represents the subprocess $\defl_2 \circ\sh^k(\X_a)$, and $a^\internal_k$ represents $\defl_2 \circ \sh^k (\eta_a)$.}
    \label{fig: construction 2 stationary process graph}
\end{figure}

\begin{remark}
    If $G^\internal$ is derived from an underlying time series graph, then it satisfies the following structural condition: Suppose $(v^\internal, i) \to (v,j)$ and $(v^\internal, k) \to (v, i)$ are links on $G^\internal$, then $(v^\internal , k) \to (v,j)$ must be a link on $G^\internal$ as well. This follows as $( v^\internal, i) \to (v, j)$ implies that there is a directed path $\pi_{i,j} \in \paths_{\smry(\G_v)}((v,i), (v,j))$, and $(v^\internal, k) \to (v, i)$ implies a directed path $\pi_{k,i} \in \paths_{\smry(\G_v)}((v,k), (v,i))$. Hence, the concatenation of $\pi_{k,i}$ and $\pi_{i,j}$ gives a directed path $\pi_{k,j}$ from $(v,k)$ to $(v,j)$, and therefore, the link $(v^\internal, k) \to (v, j)$. In a similar way one can argue that whenever $(u,i) \to (v, j)$ and $(v^\internal, k) \to (v,l)$ on $G^\internal \cup G$, then the process graph must also contain the link $(u,i) \to (v,l)$. In particular, the graph in Figure \ref{fig: 2-stationary sep general case} cannot arise as the process graph of an underlying $2$-stationary time series graph, due to $(u,0) \not\to (v,1)$. 
\end{remark}

We define the \emph{internal dynamics} of $v \in V$ to be the $P$-stationary process $\X^\internal_v$ that satisfies 
\begin{align}\label{eq: internal dynamics P-stationary SVAR}
    \X_v^\internal = \Phi^\top_{v,v} \ast \X^\internal_v + \eta_v.
\end{align}  
The filter $\iota_v$ generating the internal dynamics (meaning that $\X^\internal_v = \iota_v^\top \ast  \eta_v$) is implicitly defined by the following linear relation: 
\begin{align}\label{eq: filter internal dynamics P-stationary SVAR}
    \iota_v &= \mathbb{I}_P + \iota_v \ast \Phi_{v,v} \iff \iota_v \ast (\mathbb{I}_P - \Phi_{v,v}) = \mathbb{I}_P. 
\end{align}
In fact, the solution to the above equation is given by $\iota_v = S(\Phi_{v,v})$. In order to make sure that $\iota_v$ is a filter, i.e., absolutely summable, we need to impose a condition on $\Phi_{v,v}$. Its explicit formulation generalizes the stability notion used in \cite{reiter2025sep, reiter2026causalII} and will be given at the end of this subsection. 

\begin{example}\label{example: internal filter}
    We revisit the example associated with the time series graph of order one shown in Figure \ref{fig: 2 stationary ts graph}. Equation (\ref{eq: filter internal dynamics P-stationary SVAR}) gives the following recursive relations for the internal dynamics of any order one $2$-stationary SVAR  process: 
    \begin{align*}
        \iota (0) &= (I_2 - \Phi(0))^{-1}, & \iota(k) &= \iota(0)(\Phi(1)\iota(0))^k. 
    \end{align*}
    For the process $u$, with $\alpha_0 \coloneqq \phi_{u,u}^{(1,0)}(0)$ and $\alpha_1 \coloneqq \phi_{u,u}^{(0,1)}(1)$, we get  
    \begin{align*}
        \iota_u(0)&= \begin{bmatrix}
            1 & 0  \\ \alpha_0 & 1
        \end{bmatrix}, &
        \iota_u(k)&= \begin{bmatrix}
            \alpha_0^{k}\alpha_1^{k}  & \alpha_0^{k-1}\alpha_1^{k} \\
            \alpha_0^{k+1} \alpha_1^{k} & \alpha_0^k \alpha_1^k 
        \end{bmatrix}. 
    \end{align*}
    For the process $v$, with $\gamma_{0,1} = \phi^{(0,1)}_{v,v}(1)$ and $\gamma_{0,0} = \phi_{v,v}^{(0,0)}(1)$ and $\gamma_{1,1} = \phi_{v,v}^{(1,1)}(1)$, we get 
    \begin{align*}
        \iota_v(0)&=\begin{bmatrix}
            1 & 0 \\ 0 &  1
        \end{bmatrix}, & \iota_v(k)&= \begin{bmatrix}
            \gamma_{0,0}^k & \gamma_{0,1}\sum_{j=1}^{k-1} \gamma_{0,0}^{j} \gamma_{1,1}^{k-1-j} \\ 0 & \gamma_{1,1}^k
        \end{bmatrix}.
    \end{align*}
\end{example}

The following states that the filters $\iota_v$, defined in Equation \ref{eq: internal dynamics P-stationary SVAR} are compatible with the graph $G^\internal$. 
\begin{proposition}\label{prop: zero filters internal dynamics SVAR}
    Let $v \in V$ and $i,j\in [0,P-1]$ be such that $(v^\internal, i) \not\to (v,j)$, then $\iota^{(i,j)}_v = 0$.
\end{proposition}

It remains to construct the \emph{link filters} compatible with the process graph $G=(V \times [0, P-1], D)$. For two components $u, v \in V$ we define the $P$-stationary filter 
\begin{align}\label{P-stationary link filter}
    \lambda_{u,v} &\coloneqq \begin{cases}
        \Phi_{u,v} \ast \iota_v & u \neq v, \\
        0 & u = v. 
    \end{cases}  
\end{align} 
Compatibility of the link filters with the process graph is a direct consequence of Proposition \ref{prop: zero filters internal dynamics SVAR} and the definition of the process graph.
\begin{corollary}
    The $P$-stationary filter $\Lambda = [\lambda_{u,v}]_{u,v \in V}$ is compatible with the process graph $G$. 
\end{corollary}
\begin{example}\label{example: link filter}
    For the $2$-stationary SVAR model with underlying time series graph shown in Figure \ref{fig: 2 stationary ts graph}, we set $\beta_0 \coloneqq \phi_{u,v}^{(0,0)}(0)$ and $\beta_1 \coloneqq \phi_{u,v}^{(0,1)}(1)$. Then the 2-stationary link filter for $u \to v$ is 
    \begin{align*}
        \lambda_{u,v}(k)&= \Phi_{u,v}(0)\iota_v(k) + \Phi_{u,v}(1)\iota_v(k-1) \\
        &= \begin{bmatrix}
            \beta_0 \gamma_{0,0}^k & \beta_0 \gamma_{0,1}\sum_{j=0}^{k-1} \gamma_{0,0}^{j} \gamma_{1,1}^{k-1-j} + \beta_1\gamma_{1,1}^{k-1} \\
            0 & 0
        \end{bmatrix}.
    \end{align*}
\end{example} 
Note that if $G$ is acyclic, and $\Phi$ is such that every $\iota_v$, as in Equation (\ref{eq: filter internal dynamics P-stationary SVAR}), is absolutely summable, then $S(\Lambda)$ is absolutely summable. However, if $G$ is cyclic, then each $\iota_v$ being absolutely summable does not suffice to ensure absolute summability of $S(\Lambda)$. In the following, we provide a set of filters $\Phi$ that yield a stable SEP representation, even if $G$ is cyclic. Let $\tilde V_1, \dots, \tilde V_K$ be the strongly connected components of the summary graph $\smry(\G)$, then we define the following subset of $\G$-compatible filters:
\begin{align}\label{stable filters}
    \R^\G_\stable & \coloneqq \left\{ \Phi \in \R^\G\mid \left\lVert [\Phi^\top]_{\tilde V_i,\tilde V_i }\right\rVert_1 < 1, \text{ for all } i \in[K] \right\} ,  
\end{align}
with $\lVert\Phi\rVert_1 \coloneqq \lVert \sum_{k\in \Z}\lvert\Phi(k)\rvert \rVert_1$, where on the right, absolute values are taken entry-wise and the norm is the operator norm w.r.t. $\lVert \bullet \rVert_1$. The following states that the solution of the SVAR model coincides with the solution of the associated SEP. Its proof is deferred to the supplementary material.  
\begin{proposition}\label{prop: P-stationary SEP representation of SVAR process}
    Let $(\G, \Phi, \eta)$ be a $P$-stationary SVAR model and $\Phi \in \R^\G_\stable$, then the associated $P$-stationary SEP $(G, \Lambda, \iota^\top \ast \eta)$ is stable, and furthermore, a $P$-stationary process $\X$ is a solution to the SVAR model if and only if it is a solution to the $P$-stationary SEP model. 
\end{proposition}
\subsection{The $RP$-stationary SEP representation of a $P$-stationary SVAR model}
Given a stable $P$-stationary SEP $\mathbb{S}=(G, \Lambda, \X^\internal)$, the properties of the embedding $\mathsf{E}_{R}$ allow one to identify an $RP$-stationary SEP $\mathbb{S}_R =(G_R , \ebdg_R(\Lambda), \ebdg_R(\X^\internal))$ that is equivalent to $\mathbb{S}$, in the sense that their solutions coincide. If the SEP $\mathbb{S}$ has an underlying $P$-stationary SVAR model $\mathcal{V} =(\G, \Phi, \eta)$, then it is possible to construct construct an $RP$-stationary time series graph $\G_R$ such that the embedded SVAR model $\mathcal{V}_R = (\G_R, \ebdg_R(\Phi), \ebdg_R(\eta))$ and $\mathcal{V}$ are equivalent to each other. However, we will not need the embedded time series graph $\G_R$ in the present work. Instead, we will make use of the following: Since the embedding preserves the identity filter and is compatible with convolution, see Proposition \ref{prop: embedding properties}, it holds that    
\begin{align}\label{eq: embedded SVAR filters}
    \mathsf{E}_{R}(\iota_v  ) &= \mathbb{I}_{R} + \mathsf{E}_{R}(\Phi_v)\ast\mathsf{E}_{R}(\iota_v),  & \mathsf{E}_{R}(\lambda_{u,v}) &= \mathsf{E}_{R}(\Phi_{u,v}) \ast \mathsf{E}_{R}(\iota_v). 
\end{align} 

\section{The $P$-stationary SEP in the frequency domain}\label{section: frequency domain}
To construct the frequency domain representation of a $P$-stationary SEP, we briefly recall the Fourier transformation, see, e.g., \cite{brockwell2009timeseries}. Then we will see that the frequency domain representation of an $P$-stationary SVAR process is parameterized over the field of rational functions with real coefficients. In the final part of this section, we will derive a linear relationship between the $P$-stationary spectrum, i.e., the Fourier transformed ACS, and its $RP$-stationary embedding. In the next section, we will use this relationship to decide whether a given entry in the $P$-stationary spectrum can be linearly recovered from a temporally incomplete observational pattern. 

Throughout this section, $(\G, \Phi, \eta)$ specifies a $P$-stationary SVAR model over a finite set of vertices $V$. Furthermore, we will assume throughout this section that $\Phi\in \R^\G_\stable$. If not stated otherwise, we will refer by $\mathbb{S}=(G, \Lambda, \iota^\top \ast \eta)$ to the associated $P$-stationary SEP representation. 
\subsection{Fourier transformation}
Let $S^1$ be the set of all $z\in \C$ such that $|z|=1$. The linear and continuous Fourier transform is a Hilbert isomoprphism between the space of complex-valued square summable sequences, denoted by $\ell^2(\Z; \C)$, and the space $L^2(S^1; \C)$ consisting of the (equivalence classes of) square integrable functions $f:S^1 \to \C$.
The Fourier transform takes a sequence $\lambda$ to the function $\F(\lambda) = \left( z \mapsto \sum_{k \in \Z} \lambda(k)z^k \right)$. 

It will be useful to recall that the Fourier transform is compatible with the (tilted) convolution of sequences in the following sense: Let $ \lambda , \mu \in \ell^2(\Z; \C)$ be sequences, then it holds that $\F(\lambda \ast \mu ) = \F(\lambda) \F(\mu)$ and $\F(\lambda \ast' \mu ) = \F(\lambda) \F(\mu)^\ast$, where right hand sides are point-wise products of functions, and $\ast$ in superscript means taking complex conjugate point-wise, see \cite{reiter2025sep}. 

For $V$ a finite set, we denote by $I_V$ the identity mapping in $L^2(S^1; \C^{V \times V})$, i.e., $I_V(z) = I_V$ is the identity matrix for every $z \in S^1$. We omit the subscript if there is no ambiguity about the index set $V$ of the underlying vector space.  

The frequency domain analogues of the $j$-shift and $P$-inflation operators (for sequences) map a function $f \in L^1(S^1; \C)$ to $\sh^{j}(f)\coloneqq \left( z \mapsto z^{j} f(z) \right)$ and $\infl_P(f) \coloneqq \left( z \mapsto f\left(z^P\right) \right)$, respectively. 
These mappings are linear and related to the their analogues in the space of sequences as follows:
\begin{align}\label{Fourier transformed shift and stretch operator}
    \F(\sh^j(\lambda))&= \sh^j(\F(\lambda)), & \F (\infl_P(\lambda)) &= \infl_P (\F(\lambda)). 
\end{align} 
The embedding operators $\mathsf{e}_{R}, \mathsf{e}_R'$ for sequences, defined in Section \ref{section: P-stationary processes}, can be translated to the frequency domain. Specifically, since $\mathcal{F}$ is an isomorphism, we can define for $f \in L^2(S^1, \C^{V \times U} )$ and some positive integer $N$ the following: 
\begin{align}\label{def: embedding frequency domain}
    \mathsf{e}_{N}(f) &\coloneqq \mathcal{F} \circ \mathsf{e}_{N} \circ \mathcal{F}^{-1}(f), & \mathsf{e}_{N}'(f) &\coloneqq \F \circ \mathsf{e}_{N}' \circ \F^{-1}(f).  
\end{align}
Note that $\mathsf{e}_N(f), \mathsf{e}_N'(f) \in L^2(S^1, \bmat(\C, [0,P-1],(V,U)))$, i.e., they are block-structured. Furthermore, for a block-structured function $F: L^2(S^1;\bmat(\C, V , [0,P-1]))$, we define the block-wise operators $\ebdg_N(F) = [\mathsf{e}_N(F_{v,w})]_{v,w \in V}$ and $\ebdg_N'(F) = [\mathsf{e}'_N(F_{v,w})]_{v,w \in V}$. 

To derive the main result of this section, we will use that the embedding has the following properties. 
\begin{lemma}\label{lemma: properties frequency embedding}
    Let $f,g \in L^2(S^1, \C^{V \times V})$, then it holds that $\mathsf{e}_N(fg) = \mathsf{e}_N(f)\mathsf{e}_N(g)$. Furthermore, $\mathsf{e}_N(I) = \mathsf{e}_N'(I) = I$, i.e., the operators preserve the identity function. In particular, if $f$ is point-wise invertible, then it holds that $\mathsf{e}_N(f^{-1}) = \mathsf{e}_N(f)^{-1}$. Setting $\tilde f = \mathsf{e}_N(f)$, it holds that $\tilde f_{i,j} = \tilde f_{i+1,j+1}$ for all $i,j \in [0,N-2]$, and $\tilde f_{N-1, k-1} = \sh^1(\tilde f_{0,k})$, for every $k \in [0,N-1]$. Finally, $\mathsf{e}_N'(f) = [\tilde f_{j,i}]_{i,j \in [0,N-1]}$, where the right hand side is the block-transpose of $\tilde f$.  
\end{lemma}
The statements in Lemma \ref{lemma: properties frequency embedding} follow directly from the definition of the embedding in the sequence space and Lemma \ref{lemma: (tilted) convolution under the embedding} and the basic properties of the Fourier transform.

Suppose $\X$ is a $P$-stationary process and let $\Sigma$ be its $P$-stationary ACS. If $\Sigma$ is absolutely summable, then we define the spectral density of $\X$ to be $\fcum \coloneqq  \F(\Sigma)$.

\subsection{$P$-stationary SVAR frequency domain representation}
Let $\mathbb{S}=(G, \Lambda, \iota^\top \ast \eta)$ be a stable $P$-stationary SEP, then we define $H \coloneqq \F(\Lambda)$ to be the block-matrix of link functions, and $f = \F(\iota)$ to be the block-matrix of internal functions. Equation (\ref{equation: acs p-stationary sep}) in combination with the basic properties of the Fourier transform allow to factorize the spectrum of the solution process $\X$ and the internal dynamics $\X^\internal = \iota^\top \ast \eta$, respectively, as follows: 
\begin{align*}
    \fcum &= (I - H)^{-\top} \fcum^\internal (I- H)^{-\ast}, & \fcum^\internal &= \fcum^\internal = f^\top \Omega f^\ast,
\end{align*}
where $\Omega$ denotes the constant function that maps every $z\in S^1$ to the diagonal matrix $ \acs(\eta)(0)$.
In the following, $(G, H, \fcum^\internal)$ will be referred to as the \emph{spectral representation} of $\mathbb{S}$.

We will now see that the spectral representation of the $P$-stationary SEP associated with the SVAR model $(\G, \Phi, \eta)$ can be formulated over the rational functions. First, recall that the space of polynomials with real coefficients in one indeterminate $z$ is 
\begin{align*}
    \R[z] & \coloneqq \left\{ \sum_{k=0}^n \alpha_k z^k \mid \alpha_k \in \R, \alpha_n \neq 0 \right\}. 
\end{align*}
The field of fractions of $\R[z]$ is the field of rational functions and will be denoted $\R(z)$, on which there is an involution $(\bullet)^\ast:\R(z) \to \R(z)$, i.e., a self-inverse isomorphism that respects addition and multiplication of rational functions, with the additional property that $r^\ast(z) = r(z^\ast)$, for every $z\in S^1$, see \cite{reiter2026causalII}. 
For $v \in V$, it follows from the definition of the internal filter (Equation \ref{eq: internal dynamics P-stationary SVAR}) and the compatibility of the Fourier transform with convolution that 
\begin{align}\label{P-stationary internal function}
    f_v(\Phi) &= (I_P - \F(\Phi_{v,v}))^{-1} \in \R(z)^{P \times P}. 
\end{align}
This matrix has entries in $\R(z)$ since $\F(\Phi_{v,v})$ is a matrix of polynomials. Combining Equation (\ref{P-stationary internal function}) with the definition of the link filter, i.e., Equation (\ref{P-stationary link filter}), yields that the matrix-valued link functions associated with two distinct processes $u \neq v$ is given by the matrix of rational functions
\begin{align}\label{P-stationary link function}
    h_{u,v}(\Phi) &= \F(\Phi_{u,v})f_v = \F(\Phi_{u,v})(I_P - \F(\Phi_{v,v}))^{-1} \in \R(z)^{P \times P}. 
\end{align}
\begin{example}
    For the time series graph from Figure \ref{fig: 2 stationary ts graph}, and notation as in Example \ref{example: internal filter} and Example \ref{example: link filter}
    the internal function of the component $v$ takes the following form 
    \begin{align*}
        f_v(z) &= \frac{1}{(1-\gamma_{0,0}z)(1-\gamma_{1,1}z)} \begin{bmatrix}
            1- \gamma_{1,1}z & \gamma_{0,1} \\
            0 & 1- \gamma_{0,0} z
        \end{bmatrix},
    \end{align*}
    and the 2-stationary link function of $u \to  v $ is 
    \begin{align*}
        h_{u,v}(z) &= \begin{bmatrix}
            \beta_0  & \beta_1z \\ 0 & 0
        \end{bmatrix} f_v(z) = \frac{1}{(1- \gamma_{0,0}z)(1-\gamma_{1,1}z)}\begin{bmatrix}
            \beta_0(1- \gamma_{1,1}z) & \beta_0 \gamma_{0,1}  + \beta_1(z-\gamma_{0,0}z^2) \\ 0 & 0
        \end{bmatrix} .
    \end{align*}
\end{example} 
For the following, we define $\K_\G \coloneqq \R(\Phi, \Omega, z)$ to be the field of rational functions with indeterminates $\Phi, \Omega$ and $z$, where $\Phi$ is the vector of indeterminates where each entry corresponds to a unique link on $\G$. Similarly, we treat $\Omega$ as a vector of indeterminates.
\begin{remark}\label{remark: parameterized spectral parameterization} 
In the next section, it will be useful to consider the frequency domain objects associated with a time series graph $\G$ as parameterized rational functions, i.e., 
\begin{align*}
    H, f, \fcum^\internal, \fcum \in \bmat(\K_\G, V, [0,P-1]).     
\end{align*}
\end{remark}

\subsection{$RP$-stationary frequency representation of a $P$-stationary SVAR process}
Let $(\G, \Phi, \eta)$ be an $P$-stationary SVAR model. The goal of this subsection is to show that the associated $P$-stationary spectrum $\fcum$ can be linearly recovered (over $\R(z)$) from its $RP$-stationary embedding $\ebdg_{R}(\fcum)$. This linear relationship will be central in the derivation of the identifiability criterion for link functions in $P$-stationary SVAR processes under temporally incomplete observations. 

To capture the general structure of the embedded spectrum $\ebdg_R(\fcum)$, we define the following matrices over $\R(z)$. For a non-zero $r \in \R(z)$ and $P\times P$-dimensional matrices $A_0, \dots , A_{R-1}$ with entries in $\R(z)$, we define the block-matrix 
\begin{align}\label{eq: special Toeplitz matrices}
    \mathbf{A}(A_0, \dots , A_{R-1}) &\coloneqq \begin{bmatrix}
        \mathbf{A}_{i,j}
    \end{bmatrix}_{i,j \in [0,R-1]}, & \mathbf{A}_{i,j} & \coloneqq \begin{cases}
        A_{j -i } & j \geq i, \\
        rA_{P-i +j} & j < i. 
    \end{cases}
\end{align}
We will denote the set of all block-matrices that possess the structure (\ref{eq: special Toeplitz matrices}) by $\mathcal{T}_r(R, P) \subset \bmat(\R(z), [0,R-1], [0,P-1])$. The following proposition states some useful properties of $\mathcal{T}_r(R,P)$. We provide the proof in the supplementary material.
\begin{proposition}\label{prop: field of special Toeplitz matrices}
    The set $\mathcal{T}_r(R,P)$ is multiplicative and additive, i.e., for $\mathbf{B}, \mathbf{C} \in \mathcal{T}_r(R, P)$, it holds that $\mathbf{B}+ \mathbf{C},\mathbf{B} \mathbf{C} \in \mathcal{T}_r(R,P)$. Furthermore, if $\mathbf{B}$ is invertible, then $\mathbf{B}^{-1}\in \mathcal{T}_r(R,P)$. 
\end{proposition}
For $A = [A_{i,j}]_{i,j \in [0,R-1]} \in \mathcal{T}(R,P)$, we define  $A^{\top_R} = [A_{j,i}]_{i,j \in [0,R-1]}$ to be the block-transpose, i.e., the arrangement of the blocks is transposed while the arrangement of the entries within the blocks is left unchanged.   
\begin{lemma}\label{lemma: transposition and conjugation isomorphism}
    For a positive integer $P$, and a non-zero $r\in \R(z)$, there are the following bijections
    \begin{align*}
        (\bullet)^{\top_R} &: \mathcal{T}_r(R,P) \longrightarrow \mathcal{T}_{r^{-1}}(R,P), \quad A \mapsto A^{\top_R},&
        (\bullet)^\ast & : \mathcal{T}_r(R,P) \longrightarrow\ \mathcal{T}_{r^\ast}(R,P), \quad A \mapsto A^\ast,
    \end{align*}
    where transposition induces an anti-isomorphism and conjugation induces an isomorphism (w.r.t. to matrix multiplication). 
\end{lemma}
For $k \in \Z$, we will abbreviate $\mathcal{T}_{(z \mapsto z^k)}(R,P)$ by $\mathcal{T}_{k}(R,P)$. In the following, we will call a matrix $F \in \bmat(\R(z), V, [0,RP-1])$ block-wise $\mathcal{T}_r(R,P)$ if for every $u,v\in V$, it holds that $F_{u,v} \in \mathcal{T}_r (R,P)$. For the following proposition, note that the identities in (\ref{eq: embedded SVAR filters}) imply that the embedded link functions $\ebdg_R(H )$ and internal functions $\ebdg_R(f)$ are matrices over $\R(z)$. 
\begin{proposition}\label{corollary: inflated structures in the frequency domain}
     The embedded matrix of link function $\ebdg_{R}(H)$, the embedded internal spectrum $\ebdg_{R}(\fcum^\internal)$, and the embedded spectrum $\ebdg_R(\fcum)$ are block-wise $\mathcal{T}_1(R,P)$.  
\end{proposition}
\begin{proof}
    That $\ebdg_R(H)$ is block-wise $\mathcal{T}_1(R, P)$ is a direct consequence of Lemma \ref{lemma: properties frequency embedding}. That the embedded internal spectrum is block diagonal follows from the structure of the filter $\iota$ and the fact that the components of the white noise process $\eta$ are mutually independent. To see that each block on the diagonal is in $\mathcal{T}_1(R, P)$, fix some $v \in V$. First, we apply Proposition \ref{prop: embedding properties} to Equation (\ref{equation: acs p-stationary sep}) and then apply the Fourier transform and use its basic properties mentioned at the beginning of this section to get the following factorization:
    \begin{align*}
        \mathsf{e}_R (\fcum^\internal_v) &= \mathsf{e}_R(f_v^\top) \mathsf{e}_R(\Omega) \mathsf{e}'_R(f_v)^\ast.
    \end{align*}
    It follows from Lemma \ref{lemma: properties frequency embedding} and Lemma \ref{lemma: transposition and conjugation isomorphism} that every factor in this factorization is an element in $\mathcal{T}_1(R,P)$. Hence, it follows from Proposition \ref{prop: field of special Toeplitz matrices} that $\ebdg_R(\fcum^\internal)$ is block-wise $\mathcal{T}_1(R,P)$. 

    It remains to show that $\ebdg_R(\fcum)$ is block-wise $\mathcal{T}_1(R,P)$. For the following, we set $\tilde H \coloneqq  (I - H )^{-1}$.
    Since the embedding $\ebdg_R(\bullet)$ commutes with taking inverse, it follows from what we have just shown, and the Cayley-Hamilton theorem combined with Proposition \ref{prop: field of special Toeplitz matrices} that $\ebdg_R(\tilde H^\top )$ and $\ebdg_R(\tilde H^\ast)$ are block-wise $\mathcal{T}_1(R,P)$.
    To see that the embedded spectrum $ \ebdg_{R}(\fcum)$ is block-wise $\mathcal{T}_1(R,P)$, we factorize 
    \begin{align}\label{equation: embedded spectral factorization}
        \ebdg_{R}(\fcum) &= \ebdg_{R}(\tilde H^{\top}) \ebdg_{R}(S^\internal) \ebdg_{R}'(\tilde H^{\ast}).
    \end{align}
    Since each factor is block-wise $\mathcal{T}_1(R,P)$, it follows from Proposition \ref{prop: field of special Toeplitz matrices} that the product, i.e., $\ebdg_R(\fcum)$ is block-wise $\mathcal{T}_1(R,P)$.
\end{proof}
\begin{example}
    Considering the same setting and notation as in Example \ref{example: internal filter}, the embedded internal functions for the process $u$ and $v$, respectively, are
    \begin{align*}
        \mathsf{e}_2(f_u)(z)&= \begin{bmatrix}
            1 & -\alpha_0 & 0 & 0  \\
            0 & 1 & - \alpha_1 & 0 \\
            0 & 0 & 1 & -\alpha_0 \\
            -\alpha_1 z &0 & 0 & 1
        \end{bmatrix}^{-1}, & \mathsf{e}_2(f_v)(z) &= \begin{bmatrix}
            1 & 0 & - \gamma_{0,0} & -\gamma_{0,1} \\ 
            0 &  1 & 0 & - \gamma_{1,1}  \\ 
            -\gamma_{0,0} z & -\gamma_{0,1}z & 1 & 0 \\ 
            0 & - \gamma_{1,1} z & 0 & 1 
        \end{bmatrix}^{-1}.
    \end{align*}
\end{example} 

We will now establish a linear relationship between the $P$-stationary spectrum and its $RP$-stationary embedding. 
First, recall that the shift and inflation operators $\sh^j$ and $\infl_R$ map rational functions to rational functions. Hence, entry-wise application yield the following homomorphisms of rings, i.e., they are compatible with matrix addition and multiplication: 
\begin{align*}
    \sh^j: \mathcal{T}_k (R,P) \to \mathcal{T}_{k}(R,P), \quad \infl_R: \mathcal{T}_k(R,P) \to \mathcal{T}_{kR}(R,P). 
\end{align*}  
To state the linear relationship in a way that will be useful in the next section, we define for a given set $I \subset [0,R-1]^2$ the associated mapping 
\begin{align}\label{mapping: spectral identifiability}
    \varphi_I : I \to [0,R-1], \quad (i,j) \mapsto (j-i) \mathbf{1}_{i \leq j}(i,j) + (R-i +j)\mathbf{1}_{i > j}(i,j).
\end{align}
\begin{proposition}\label{prop: deflation in the frequency domain}
    Let $u,v \in V$ and $I \subset [0,R-1]^2$ be such that the mapping $\varphi_I$ is a bijection, then for $\tilde \fcum = \mathsf{e}_R (\fcum_{u,v})$, it holds that 
    \begin{align*}
        \fcum_{u,v} &= \sum_{(i,j) \in I}\sh^{j-i} \circ \infl_R\left(\tilde \fcum_{i, j}\right) .
    \end{align*}
\end{proposition}
\begin{proof}
First, it follows from Lemma \ref{lemma: linear recovery of the embedded filter} and the two identities in (\ref{Fourier transformed shift and stretch operator}) that
\begin{align}\label{equation: spectral recovery}
    \fcum_{u,v} &= \sum_{j=0}^{R-1} \sh^j \circ \infl_R (\tilde \fcum_{0,j}).
\end{align}
Secondly, let us define $I_\leq = \{ (i,j) \in I \mid i \leq j \} $ and $I_> = \{ (i,j) \in I \mid i > j \}$. Note that it follows from Lemma \ref{lemma: properties frequency embedding} that for $(i,j) \in I_\leq$, it holds that $\tilde \fcum_{i,j} = \tilde \fcum_{0, j-i}$ and for $(i,j) \in I_>$, it holds that $\tilde \fcum_{i,j} = \sh(\tilde\fcum_{0,R-i+j})$. Using the assumption on $I$, the identity $\infl_R \circ \sh^j = \sh^{Rj} \circ \infl_R$ and the decomposition (\ref{equation: spectral recovery}), we compute  
\begin{align*}
    \fcum_{u,v} &= \sum_{(i,j) \in I_\leq } \sh^{j-i} \circ \infl_R (\tilde \fcum_{0,j -i}) + \sum_{(i,j) \in I_>} \sh^{R-i +j} \circ \infl_R (\tilde \fcum_{0, R -i +j }) \\
    &=\sum_{(i,j)\in I_\leq} \sh^{j-i} \circ \infl_R (\tilde \fcum_{i,j}) + \sum_{(i,j) \in I_>} \underbrace{\sh^{R-i+j} \circ \infl_R (\sh^{-1}(\tilde \fcum_{i,j}))}_{= \sh^{j-i} (\tilde\fcum_{i,j})} \\
    &= \sum_{(i,j) \in I} \sh^{j-i} (\tilde \fcum_{i,j}).
\end{align*}
\end{proof}
\section{Identification under temporally incomplete observations}\label{section: identification}
Throughout this section, $(\G, \Phi, \eta)$ denotes a $P$-stationary SVAR process whose components are indexed by a finite set $V$, and $\Phi \in \R^\G_\stable$. First, we make precise what we mean by temporally incomplete. Then we will derive a simple criterion to check whether a given entry in the $P$-stationary spectrum can be recovered from a given temporally incomplete observation. Furthermore, we provide a graphical criterion to determine the invertibility of a given subspectrum, which is in fashion of \cite{sullivant2010trek, reiter2026causalII}, and which we will use to derive a instrumental process based graphical criterion to decide the identifiability of a set of link functions under temporally incomplete observations. 
\subsection{Semi-periodic subsets of $\Z$}
\begin{definition}[Semi-periodic set]
    Let $P$ be a non-negative integer. A subset $O \subset \Z$ is said to be \emph{$P$-periodic} if there is some $m \in [0,P-1]$ such that $O =P\Z-m$. A set of integers $O$ is called \emph{periodic} if there exists an integer $P$ such that $O$ is $P$-periodic. 
    A \emph{semi-periodic set} is a finite union of periodic sets. A semi-periodic set $O$ is said to have \emph{common periodicity} $P$ if it can be expressed as a union of $P$-periodic sets. 
\end{definition} 
Given positive integers $P_1, \dots, P_d$, recall that their \emph{least common multiple} is given by the smallest positive integer $k$ such that $k \in P_i \Z$ for every $i \in [d]$. 
\begin{lemma}
    Let $O = \bigcup_{i =1}^d P_i \Z -m_i$ be a semi-periodic set, where $P_i >0$ and $m_i \in [0,P_i-1]$ for all $i \in [d]$. Then the least  common multiple of the $P_i$'s is the smallest common periodicity of $O$. 
\end{lemma}
\begin{proof}
    The claim follows from the following simple observation. Suppose $O= P\Z-m$ is a $P$-periodic set, then $O$ can be written as semi-periodic set with common periodicity $RP$, where $R$ is a positive integer. Specifically, it holds that $O = \bigcup_{k = 0}^{R-1} RP \Z - m -kP$.
\end{proof}
\begin{lemma}\label{lemma: base periodicity}
    Let $\mathbf{O}=(O_v \mid v \in V)$ be a finite collection of semi-periodic sets and let $P_v $ be the smallest common periodicity of $O_v$. For $P$ the least common multiple of the positive $P_v$'s, it holds that $P$ is the smallest number such that each $O_v$ has common periodicity $P$.  
\end{lemma}
\begin{definition}
    Let $\mathbf{O}= (O_v \mid {v\in V})$ be a finite collection semi periodic sets. A positive integer $P$ is said to be a common periodicity for $\mathbf{O}$ if $P$ is a common periodicity for every $O_v$. The \emph{minimal periodicity} of $\mathbf{O}$ is the smallest common periodicity.  
\end{definition}  
In the following, we will represent finite collections of semi-periodic sets as follows. Let $\mathbf{O}= (O_v \mid {v \in V})$ be a collection of semi-periodic sets with common periodicity $P$, then define for every $v \in V$ the set $I_v \coloneqq \{ i \in [0,P-1] | P \Z-i \subseteq O_v \}$. We will refer to $\mathbf{I} =(I_v \mid v \in V)$ as the finite representation of $\mathbf{O}$ with respect to the common periodicity $P$.
\begin{example}
    Consider $V= \{u,v \}$, with semi-periodic sets 
    \begin{align*}
        O_u &= 2 \Z \cup 4\Z-1 , &O_v &= 3\Z-1 \cup 2\Z-1.
    \end{align*} The minimal common periodicity of $O_u$ is $4$ and the minimal common periodicity of $ O_v$ is 6. Consequently, the minimal common periodicity for $(O_u, O_v)$ is 12, giving the representation
    \begin{align*}
        I_u &= \{0,2,3,4,6,7,8,10,11\}, &
        I_v &=\{1,2,3,5,7,8,9, 11\}.  
    \end{align*}
\end{example}
\subsection{Spectral identifiability under temporally incomplete observations}
\begin{definition}
    Let $\X= (\X_v)_{v\in V}$ be a $V$-variate process. A collection of semi-periodic sets $\mathbf{O}=(O_v \mid v \in V)$ specifies the observational structure of the process $\X$ if $\X_v(t)$ is observable if and only if $ t \in O_v$. The observable process is $\X(\mathbf{O}) = (\X_v(O_v))_{v \in V }$. The observational structure $\mathbf{O}$ is called \emph{temporally incomplete} if there is at least one $v \in V$ such that $O_v \neq \Z$. 
\end{definition}
Let $\X$ be an $P$-stationary SVAR process compatible with $(\G, \Phi, \eta)$, and let $\mathbf{O}$ be the observational structure for $\X$. Furthermore, let $RP > 0$ be the least common multiple of $P$ and the minimal periodicity of $\mathbf{O}$. In what follows, $\mathbf{I}$ will be the finite representation of $\mathbf{O}$ w.r.t. the common periodicity $RP$.  

Let $u,v \in V$ be the indices of two component processes. To ease notation, we write $\tilde\fcum \coloneqq \mathsf{e}_R(\fcum_{u,v})$ for the embedded spectrum, which is an element in $\mathcal{T}_1(R,P)$, see Proposition \ref{corollary: inflated structures in the frequency domain}. That means, $\tilde\fcum$ is entirely specified by the submatrices $\tilde\fcum_{0,0}, \dots , \tilde \fcum_{0,R-1}$. If each of these blocks can be recovered from the \emph{observational subspectrum} $\tilde\fcum(\mathbf{I})\coloneqq [\tilde\fcum]_{I_u, I_v}$, then one can use Proposition \ref{prop: deflation in the frequency domain} to get the $P$-stationary spectral dependence $\fcum_{u,v}$. 

In view of our main result, we will need to decide wether a specific entry in the $P$-stationary spectrum, i.e., $\fcum_{u,v}^{(j_u, j_v)}$, where $j_u , j_v \in [0,P-1]$, can be recovered from the observational subspectrum $\tilde\fcum(\mathbf{I})$. For the purpose of deriving an associated criterion, we represent every $i \in [0,RP-1]$ by the tuple $(\lfloor i /P \rfloor,i \mod P)\in [0,R-1]\times [0,P-1]$. Consequently, we represent the sets $I_u, I_v$ as subsets of $[0,R-1]\times [0,P-1]$. Furthermore, for $j\in [0,P-1]$, we define the set $I_u(j) \coloneqq \{ i \in [0,R-1] \mid (i,j) \in I_u \}$ and, analogously, the set $I_v(j)$. For $j_u,j_v \in [0,P-1]$, we write $\varphi_{u,v}^{(j_u, j_v)}$ for the mapping $\varphi_{I_u(j_u) \times I(j_v)}$, as defined in (\ref{equation: spectral recovery}). Note that it follows from Proposition \ref{prop: deflation in the frequency domain} that the entry $\fcum^{(j_u, j_v)}_{u,v}$ can be recovered if $\varphi_{u,v}^{(j_u, j_v)}$ maps surjectively into $[0,R-1]$. Based on this observation, we introduce the following terminology. 
\begin{definition}\label{def: linear identifiability spectrum}
    For $u,v \in V$ and $j_u, j_v \in [0,R-1]$, the entry $\fcum_{u,v}^{(j_u, j_v)}$ is linearly identifiable under the observational structure $\mathbf{I}$ if the mapping $\varphi_{ u,v}^{(j_u, j_v)}$ is a surjection into $[0,R-1]$. The $P$-stationary spectrum $\fcum_{u,v}$ is said to be identifiable if all its entries are identifiable. 
\end{definition} 
\begin{proposition}[Identifiability conditions]\label{prop: necessary condition}
    Let $j_u, j_v \in [0,P-1]$, then $|I_u(j_u)||I_v(j_v)| \geq R$ is necessary and $|I_u(j_u)| + |I_v(j_v)| > R$ is sufficient for identifiability of the entry $\fcum_{u,v}^{(j_u, j_v)}$. 
\end{proposition}
\begin{proof}
    We begin with the necessary condition. Denoting the image of $\varphi_{u,v}^{(j_u, j_v)}$ by $\mathrm{im}_{u,v}^{(j_u, j_v)}$, the proposition follows from $$|\textrm{im}_{u,v}^{(j_u, j_v)} |\leq |I_u(j_u)||I_v(j_v)|. $$ 
    
    We continue with the sufficient condition. It suffices to consider the case $P=1$, where $I_u, I_v \subset [0,R-1]$, and we set $\varphi \coloneqq \varphi_{ u,v}$. We show the claim by contradiction. Suppose there is some $l \in [0,R-1]$ such that 
$l \notin \mathrm{im}(\varphi)= \varphi(I_u, I_v )$. This means that for every $(i,j)\in I_u \times I_v$, we have 
\begin{align*}
    l \neq \varphi(i,j) &= \begin{cases}
        j-i & i \leq j, \\
        R-i+j & i > j. 
    \end{cases}
\end{align*}
To arrive at a contradiction, we define the set $ A = \{ i+l \mid i \in I_u \cap [0, R-l]\} $ and the set $B = \{ i+l - R  \mid i \in [R-l +1, R-1] \cap I_u\} $. First, we observe that $A \cap B = \emptyset$. If there was an element in the intersection, then there would be $i \in [0,R-l]$ and $i' \in [R-l +1, R-1]$ such that 
\begin{align*}
    i + l &= i' +l- R \iff R = i' -i. 
\end{align*}
Since the equation on the right hand cannot hold, it follows that $A \cap B = \emptyset$. Furthermore, note that $ |A \cup B| = |I_v|$, and $ l \notin \mathrm{im}(\varphi)$ means that $A \cup B \cap I_v = \emptyset$, which contradicts the assumption $|I_v| + |I_u| = |I_v|+ |A \cup B | > R$. Hence,the mapping $\varphi$ must be surjective. 
\end{proof}
\begin{figure}
    \centering
    \includegraphics[width=\linewidth]{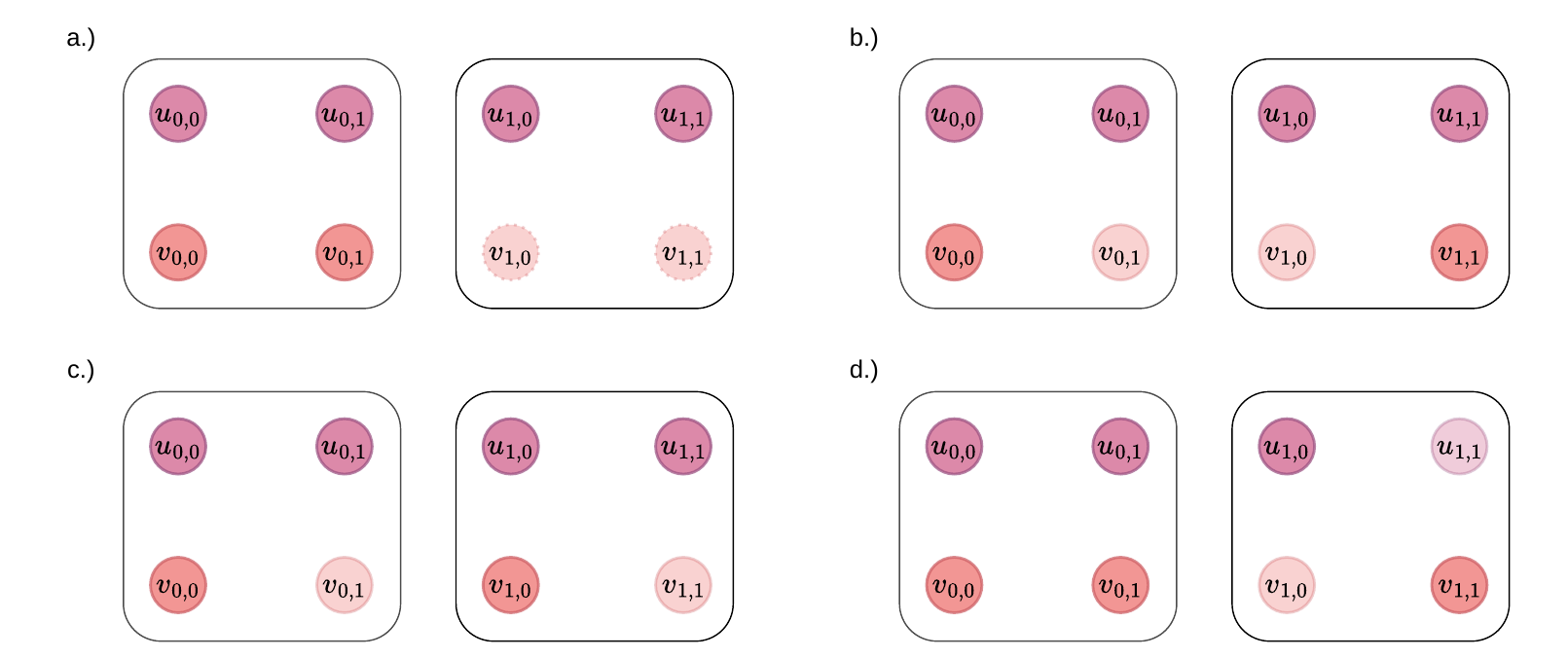}
    \caption{\small For $R=2=P$, each figure shows an observational structure of a 2-variate process, where $u_{i,j}$ represents the subprocess $\defl_4 \circ\sh^{2i+j}(\X_u)$, and analogously for $v$. The fully and palely colored nodes correspond to the observed and unobserved subprocesses, respectively.}
    \label{fig: observational structures}
\end{figure}
\begin{example}\label{example: spectral identifiability}
    In this example, we classify for $R=2$ all the partially observed scenarios where the $P$-stationary $\fcum_{u,v}$ is identifiable.
    The embedded spectrum $\tilde\fcum = \mathsf{e}_2(\fcum_{u,v})$ is block-structured as follows: 
    \begin{align*}
        \tilde\fcum &= \begin{bmatrix}
            \tilde\fcum_{0,0} & \tilde\fcum_{0,1} \\ 
            \tilde \fcum_{1,0} & \tilde\fcum_{1,1}
        \end{bmatrix}, & \tilde \fcum_{1,0} &= \sh^1 (\tilde \fcum_{0,1}). 
    \end{align*}
    Using Proposition \ref{prop: deflation in the frequency domain}, we get the following expressions for the $P$-stationary spectrum:  
    \begin{align*}
        \fcum_{u,v}(z) &= \tilde\fcum_{0,0}(z^2) + z \tilde\fcum_{0,1}(z^2) = \tilde \fcum_{0,0}(z^2) + z^{-1}\tilde \fcum_{1,0}(z^2).  
    \end{align*}
    Using Proposition \ref{prop: necessary condition}, we see that identifiability of $\fcum_{u,v}$ requires, for all $j_u, j_v \in [0,P-1]$ that $|I_u(j_u)||I_v(j_v)| \geq 2$. In particular, for every $j_u, j_v\in [0,1]$, both $I_u(j_u)$ and $I_v(j_v)$ must contain at least one element while at least one of them must contain two elements. In other words, for $R=2$, the necessary condition and the sufficient condition in Proposition \ref{prop: necessary condition} are equivalent. Examples for observational structures, with $R=2=P$, can be found in Figure \ref{fig: observational structures}, where $\fcum_{u,v}$ is identifiable under the observational structure depicted in a.) and b.), and not identifiable under the observational structures shown in c.) and d.). 

\end{example}
\subsection{Invertible subspectra} 
In this subsection, we will read all frequency domain objects associated with an $P$-stationary SVAR process as matrices over $\K_\G$, i.e., $H , f , \fcum^\internal, \fcum \in \bmat(\K_\G, V, [0,P-1])$, see Remark \ref{remark: parameterized spectral parameterization}. The goal of this subsection is to derive a graphical criterion for the summary graph $\smry(\G)$, see Definition \ref{definition: summary graph}, to decide whether a given square shaped submatrix of $\fcum$ is invertible, as a matrix over $\K_\G$. In deriving the criterion, we will exploit that the process graph imposes algebraic relations in the spectrum. To capture this relationship, we establish a \emph{trek-rule} similar to \cite{sullivant2010trek, reiter2026causalII} for $P$-stationary processes.

For a directed path $\pi = ((v_i, m_i ) \to (v_{i+1}, m_{i+1})) \in \paths_G((v,m),(w,n))$, its corresponding \emph{path-function} is defined as follows: 
\begin{align*}
    h^{(\pi)} &\coloneqq \prod_{i=1}^n h_{v_{i}, v_{i+1}}^{(m_i, m_{i+1})} \in \K_\G. 
\end{align*}
Furthermore, for a trek $\tau \in \treks_{G \cup G^\internal}((v,m), (w,n))$ with $\trektop(\tau) \in V^\internal \times [0,P-1]$, it follows from the structure of $G^\internal$ that $\tau$ must have the following form 
\begin{align*}
    (v, m) \leftarrow \dots \leftarrow (u, i) \leftarrow (u^\internal, k) \to (u,j) \to  \dots \to (w,n).   
\end{align*}
So for any trek $\tau \in \treks_{G \cup G^\internal}((v,m), ((w,n))$, we define the associated \emph{trek-function} as follows:
\begin{align}
    \fcum^{(\tau)} &\coloneqq \begin{cases}
         h^{(\pi_L)}  f_u^{(k,i)} \omega_u^{(k)} \left(f^{(k,j)}_u  h^{(\pi_R)}\right)^\ast &  \trektop(\tau) =(u^\internal,k) \in V^\internal \times [0,P-1], \\
         0 & \trektop(\tau) \in V\times [0,P-1],
    \end{cases}
\end{align}
where $\omega^{(k)}_u \in \K_\G$ is the monomial of degree one associated with the variance of $\eta_u^{(k)}$. The following describes the spectrum of a $P$-stationary SVAR process in terms of trek functions. 
\begin{proposition}[$P$-stationary trek-rule]\label{prop: trek-rule}
    Let $\G$ be a $P$-stationary time series graph, and let $G$ be the associated $P$-stationary process graph, and let $G^\internal$ the associated graph for the $P$-stationary internal dynamics. For every $(v,m), (w, n) \in V \times [0,P-1]$, it holds that 
    \begin{align*}
        \fcum^{(m,n)}_{v,w} &= \sum_{\tau \in \treks_{G \cup G^\internal }((v,m),(w,n))} \fcum^{(\tau)}.
    \end{align*}
\end{proposition}

Let $X , Y \subset V \times [0,P-1]$ be (possibly intersecting) subsets with $|X| = |Y|=N$. 
A \emph{trek system} \cite{sullivant2010trek} on the summary graph $\smry(\G)$, between $X$ and $Y$, is a tuple of treks $T=(\tau_i \mid i \in [N])$, such that there is a bijection $\alpha :[N] \to [N]$ and $\tau_i \in \treks_{\smry(\G)}(x_i, y_{\alpha(i)})$. This system $T$ is said to be without intersection if none of the $\tau_i$ intersect at any vertex. In the following, $\treks_{\smry(\G)}(X,Y)$ will denote the set of trek systems without intersection. 
\begin{proposition}\label{prop: non-intersecting paths}
    Let $\G$ be a $P$-stationary time series graph and $\smry(\G)$ its summary graph. Let $X,Y \subset V \times [0,P-1]$ be two subsets with equal cardinality, then
    \begin{align*}
        \treks_{\smry(\G)}(X,Y) \neq \emptyset \implies \det([\fcum]_{X,Y}) \neq 0.
    \end{align*} 
\end{proposition}
\begin{remark}
    Whether $[\fcum_\G]_{X,Y}$ is invertible cannot be decided by identifying non-intersecting trek systems on the process graph $G$. We demonstrate this using $\{u,v\}$-variate 2-stationary SVAR process with the following non-zero coefficients:
    \begin{align*}
        \phi_{u,v}^{(0,0)}(0)&= \beta_0, & \phi_{u,v}^{(1,0)}(1)&= \beta_1, & \phi_{v,v}^{(0,1)}(0) &= \gamma.
    \end{align*}
    Consequently, the associated summary graph is specified by the set links $\{ u_0 \to v_0, u_1 \to v_0, v_0 \to v_1\}$, which means that $\treks_{\smry(\G)}(\{u_0, u_1 \}, \{ v_0, v_1 \})= \emptyset$. On the other hand, the process graph $G$ contains the links $u_0 \to v_0$ and $u_1 \to v_1$, which means that $\treks_G(\{ u_0, u_1 \}, \{v_0, v_1 \}) \neq \emptyset$. But the link function associated with $u \to v$ has rank one, as can be seen from  
    \begin{align*}
        h_{u,v}(z) &= \begin{bmatrix}
            \beta_0 & 0 \\ \beta_1 z  & 0
        \end{bmatrix} \begin{bmatrix}
            1 & \gamma \\ 0 & 1 
        \end{bmatrix}. 
    \end{align*}
    Consequently, for $X = \{u \}\times \{0,1\} $ and $Y = \{ v \} \times \{0,1\}$, the associated subspectrum $[\fcum ]_{X,Y} = h_{u,v}^\top \Omega_u h_{u,v}^\ast$ has rank one, i.e., it is not invertible. This shows why the criterion in Proposition \ref{prop: non-intersecting paths} needs to be stated for systems of treks on $\smry(\G)$ rather than on $G$. 
\end{remark}
\subsection{Identification with generalized instrumental processes}
\begin{figure}
    \centering
    \includegraphics[width=\linewidth]{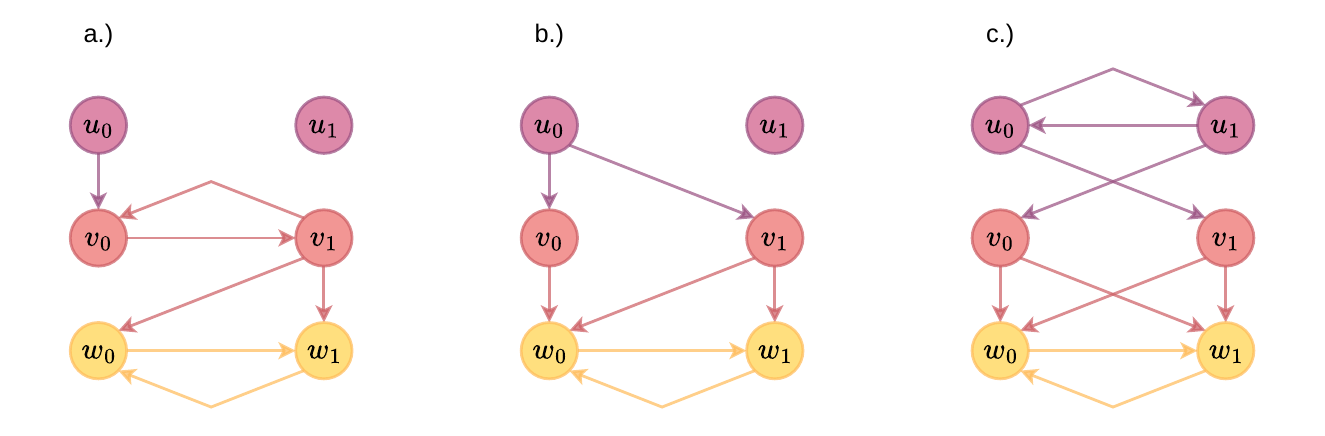}
    \caption{Three different summary graphs for the 2-stationary causal chain $u \to v \to w$. A vertex of the form $a_k$ represents the subprocess $\defl_2 \circ \sh^k(\X_a)$.}
    \label{fig: instrumental process graph causal strucutre}
\end{figure}
\begin{figure}
    \centering
    \includegraphics[width=\linewidth]{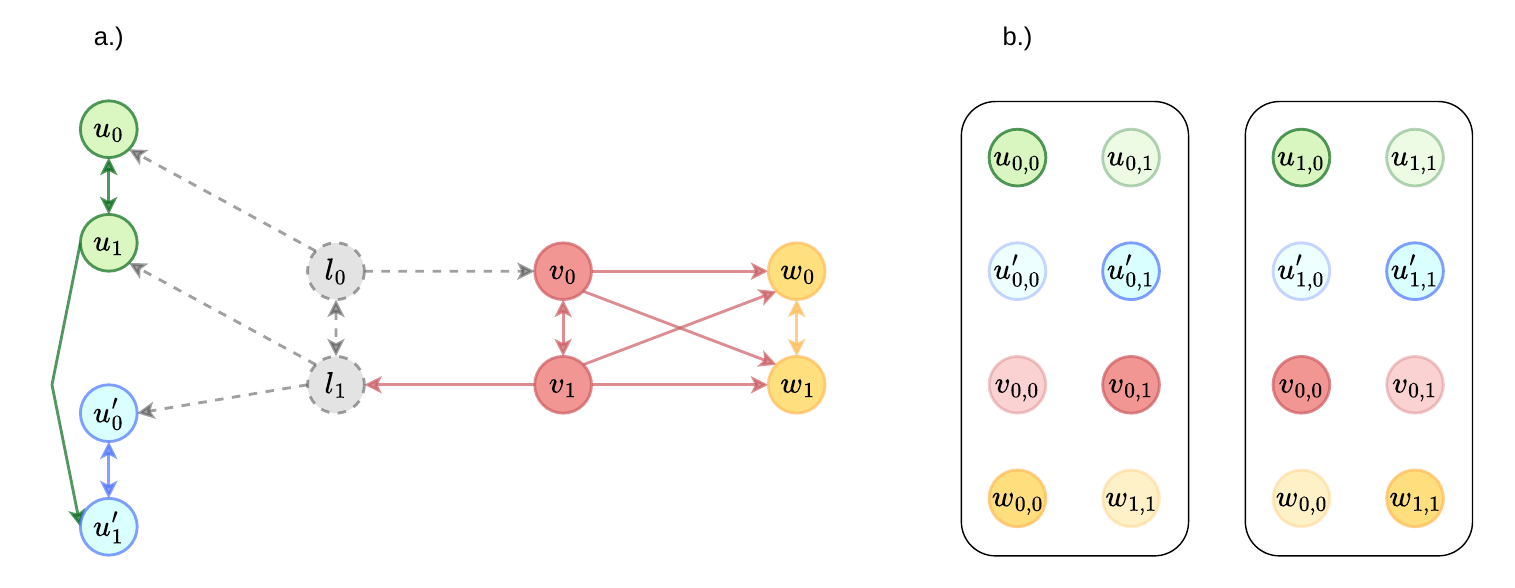}
    \caption{\small a.) The summary graph for a 2-stationary SVAR process. A vertex $a_k$ represents the process $\defl_2 \circ \sh^k(\X_a)$, and a bidirectional link represents two directed links, i.e., one in each direction. b.) The observational structure, where a vertex of the form $a_{k,j}$ represents the subprocess $\defl_4 \circ \sh^{2k +j}(\X_a)$. The fully and palely colored nodes represent the observed and unobserved subprocesses, respectively. The process $l$ is completely unobserved.}
    \label{fig: example 2 iv identification}
\end{figure}

Let $\mathbf{O} = (O_v \mid {v \in V}) \subset V \times \mathbb{Z}$ be an observational structure, and $RP$ the least common sampling rate that is a multiple of $P$, and $\mathbf{I}$ the finite representation of $\mathbf{O}$ with respect to $RP$.
\begin{question}\label{question: subsampling identifiability}
    For a process $y \in V$, does the observable spectrum $\tilde\fcum(\mathbf{I})$ provide enough information to identify the link-functions $h^{(m,n)}_{v,y}$, for all $(v,m) \in \pa_G((y,n))$?  
\end{question}
To address Question \ref{question: subsampling identifiability} algebraically, we use the following notion of identifiability, which is similar to the identifiability terms used in, e.g., \cite{foygel2012halftrek, sullivant2010trek, reiter2026causalII}. 
\begin{definition}[Rational identifiability]\label{def: rational identifiability}
    Let $\mathbf{I}$ be as above, and let $v,w \in V$ be two process indices such that $(v,i) \to (w,j)$ for some $i,j\in [0,P-1]$. The link $(v,i) \to (w,j) $ is said to be \emph{(rationally) identifiable} under $\mathbf{I}$ if there is a rational identification function 
    \begin{align*}
        \xi_{v,w}:  \K_\G^{\mathbf{I} \times \mathbf{I}} \to \K_\G,
     \end{align*} such that $\xi_{v,w}(\infl_R(\tilde\fcum(\mathbf{I}))) = h_{v,w}^{(i,j)} \in \K_\G$, where $\tilde\fcum(\mathbf{I}) = [\ebdg_R(\fcum)]_{\mathbf{I}, \mathbf{I}}$ is the observational spectrum. The $P$-stationary relation $v \to w$ is said to be identifiable if every link $(v,i) \to (w,j)$ is identifiable.
\end{definition}
\begin{figure}
    \centering
    \includegraphics[width=\linewidth]{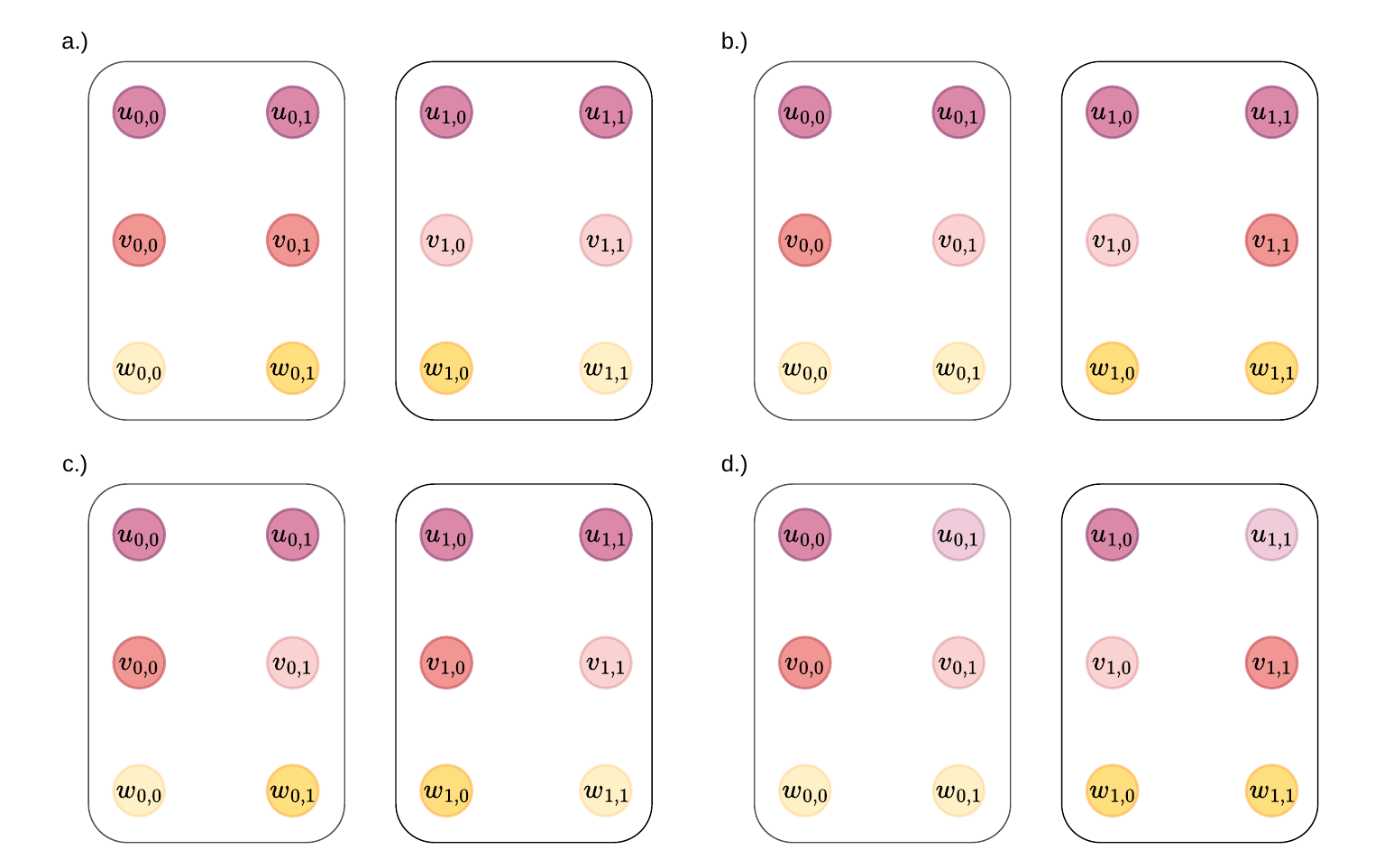}
    \caption{Four observational structures for three 2-stationary processes $u,v,w$. A vertex of the form $a_{k,j}$ represents the subprocess $\defl_2 \circ \sh^{2k+j}(\X_a)$. The fully and palely colored nodes represent the observed and unobserved subprocesses, respectively.} 
    \label{fig: instrumental identification observational structure}
\end{figure}
\begin{theorem}\label{thm: instrumental process identification}
    Let $n \in [0,P-1]$ and $y\in V$ be a process index. Suppose there is a subset $X \subset  \des_G((y,n))^c$, i.e., non-descendants of $(y,n)$ w.r.t. $G$, such that $\treks_{\smry(\G)} (X, \pa_G((y,n))) \neq \emptyset$, and such that $\fcum_{x,w}^{(i,j)}$ is identifiable under the observational structure $\mathbf{I}$, for every $(x,i) \in X $ and $(w,j) \in \pa_G((y,n)) \cup \{ (y,n)\}$. Then for every $ (w,m) \in \pa_G((y,n)) $, the link $(w,m) \to (y,n)$ is rationally identifiable under the observational structure $\mathbf{I}$.
\end{theorem}
Before proving the theorem, we briefly comment on its interpretation. Theorem \ref{thm: instrumental process identification} states that for any given time series graph $\G$, with summary graph $\smry(\G)$, a link on $G$ (recall that $G$ is determined by $\smry(\G)$) is identifiable, as of Definition \ref{def: rational identifiability}, if there is a set of instrumental (sub)-processes satisfying the conditions of Theorem \ref{thm: instrumental process identification}. Due to our notion of identifiability, it is important to note that $\R^\G_\stable$ is an open and non-empty semi-algebraic subset of $\R^\G$, meaning that it has non-zero Lebesgue measure. On the other hand, every SVAR processes (seen as a point in $\R^\G_\stable$), for which identifiability fails even though the criterion of Theorem \ref{thm: instrumental process identification} is satisfied, belongs to the set of singular points of the identification function from Definition \ref{def: rational identifiability}, and this set of singular points has Lebesgue measure zero. To summarize, if for a given  set of links, the conditions of Theorem \ref{thm: instrumental process identification} are satisfied, then for almost every stable and compatible SVAR process, the associated identification formula, applied to the observational spectrum, yields the respective link functions. 
\begin{proof}[Proof of Theorem \ref{thm: instrumental process identification}]
    Let us write $X = \{(u_i, m_i ) \mid i \in [N] \} $ and $W \coloneqq \pa_{G}((y,n)) = \{ (w_i, q_i) \mid i \in [N] \}$. Now note that every trek $\tau \in \treks_{G\cup G^\internal}((u_i, m_i), (y, n ))$ with non-zero trek function is of the form 
    \begin{align*}
        (u_i, m_i) \leftarrow  \cdots \leftarrow (z,l') \leftarrow(z^\internal, l) \to (z, l'') \to \cdots \to (w_j, q_j) \to (y,n).
    \end{align*}
    Hence, the trek rule stated in Proposition \ref{prop: trek-rule} yields the following linear relationship 
    \begin{align*}
        [\fcum]_{\{(y,n)\}, X} = [H]_{W,\{(y,k)\}}^\top [\fcum]_{W, X}.
    \end{align*}
     Using Proposition \ref{prop: non-intersecting paths} with the assumption $\treks_{\smry(\G)}(X,W)\neq \emptyset$, we conclude that 
     \begin{align*}
        [H]_{W,(y,k)}^\top &=[\fcum]_{\{(y,n)\}, X} [\fcum]_{W, X}^{-1}, 
     \end{align*}
      where by assumption all the terms on the right hand side are identifiable under $\mathbf{I}$. 
\end{proof}
\begin{example}
    We consider the scenario of a 2-stationary process consisting of three observed processes $u,v,w$, whose causal structure is encoded by the chain $u \to v\to w $. In Figure \ref{fig: instrumental process graph causal strucutre}, we show three different summary graphs encoding the chain $u \to v \to w$. In Figure \ref{fig: instrumental identification observational structure}, we show four different observational structures for the processes $u,v,w$. We will now go through these cases and check whether Theorem \ref{thm: instrumental process identification} can be applied to identify the 2-stationary effect $v \to w$. 

    a.) At the level of the process graph $G$, we have $\pa_G(w_i) = \{v_1\}$, for $i=0,1$. Under none of the observational structures, the process $v$ is completely observable, so we need to use the process $u$ as an instrumental process. Since $\paths_{\smry(\G)}(u_0, v_1) \neq \emptyset$, and since by Proposition \ref{prop: necessary condition}, the entries $\fcum^{(0,1)}_{u,v}, \fcum^{(0,0)}_{u,w}, \fcum^{(0,1)}_{u,w}$ are identifiable under all of the observational scenarios shown in Figure \ref{fig: instrumental identification observational structure}, we can apply Theorem \ref{thm: instrumental process identification} to infer that $v \to w$ is identifiable . Specifically, the identification formulas are as follows: 
    \begin{align*}
        h_{v,w}^{(1,0)} &= \frac{\fcum^{(0,0)}_{w,u}}{\fcum_{v,u}^{(1,0)}}, &
        h_{v,w}^{(1,1)} &= \frac{\fcum_{w,u}^{(1,0)}}{\fcum_{v,u}^{(1,0)}}.
    \end{align*}

    b.) In this case, we have that $\pa_G(w_i) = \{v_0, v_1 \}$, for $i= 0,1$. Hence, we need both $v_0$ and $v_1$ as instruments for the identification of $v \to w$. However, $\treks_{\smry(\G)}({u_0, u_1}\},\{v_0, v_1\}) = \emptyset$. Hence, Theorem \ref{thm: instrumental process identification} cannot be used to infer the identifiability of $v \to w$.

    c.) In this case, it also holds that $\pa_G(w_i) = \{v_0, v_1\}$, for $i = 0,1$. But now we have that $\treks_{\smry(\G)}({u_0, u_1}\},\{v_0, v_1\}) \neq \emptyset$, as $u_0 \to v_1$ and $u_1 \to v_0$ yield a non-intersecting trek system. In order to apply Theorem \ref{thm: instrumental process identification}, we need an observational structure under which $\fcum_{v,u}$ and $\fcum_{w,u}$ are identifiable. This is the case only for the observational structures shown in a.) and b.) of Figure \ref{fig: instrumental identification observational structure}. For the reasoning behind the identifiability of the spectral dependencies, the reader may review Example \ref{example: spectral identifiability}. 
\end{example}
\begin{example}\label{example: iv identification 2}
    In Figure \ref{fig: example 2 iv identification}, there are four partially observed processes $u, u', v, w$ and one latent process $l$. The 2-stationary summary graph for those  processes is depicted in Figure \ref{fig: example 2 iv identification} a.). The observational structure for is shown in Figure \ref{fig: example 2 iv identification} b.). Note that all of the processes are only partially observable. For the identification of $v \to w $, the process $u$, i.e., $u_0$ and $u_1$, could serve as an instrumental process, since $u_0 \leftarrow l_0 \to v_0$ and $u_1 \leftarrow l_1 \leftarrow v_1 $ is a non-intersecting trek system, i.e., $\treks_{\smry(\G)}(\{u_0, u_1\}, \{v_0, v_1 \} ) \neq \emptyset$. However, $u_1$ is not observable, so $\fcum_{u,v}^{(1,i)}$ is not identifiable, as of Definition \ref{def: linear identifiability spectrum}, for $i= 0,1$. Considering the summary graph, we notice that also $u_0$ and $u'_1$ can be used as an instrumental processes for the identification of $v \to w$, since $u_0 \leftarrow l_0 \to v_0 \to w_0 $ and $u_1' \leftarrow u_1 \leftarrow l_1 \to v_1 $ constitute a non-intersecting trek system. The terms $\fcum^{(0, i)}_{u, v}$ and $\fcum_{u, w}^{(0,i)}$ and $\fcum_{u', v}^{(1, i)}$ and $\fcum_{u, w}^{(1,i)}$ are identifiable for $i= 0,1$, see Proposition \ref{prop: necessary condition}. Hence, Theorem \ref{thm: instrumental process identification} implies that the 2 stationary effect $v \to w$ is identifiable. 
\end{example}
\begin{remark}
    In Example \ref{example: iv identification 2}, we used subprocesses of two different component processes to identify the 2-stationary effect $v \to w$, which was possible because we assumed a $2$-stationary process. If the process was stationary instead of 2-stationary, and if we had the same observational structure, then we could not use Theorem \ref{thm: instrumental process identification} for the identification of the then stationary effect $v \to w$, since neither $\fcum_{u, v}$ nor $\fcum_{u', v}$ are identifiable, as of Definition \ref{def: linear identifiability spectrum}. In principle, we could make use of the embedding operation, i.e., considering the equivalent 2-stationary embedding and try to use a similar algebraic structure as above (in this case, the embedded time series graph $\G_2$, mentioned earlier, will matter). However, this requires an adaptation of Proposition \ref{prop: non-intersecting paths} (since this proposition  is not stated for embedded processes), which will require additional considerations, since the embedding introduces additional relations in the link functions. This problem seems to be an interesting subject for future work. 
\end{remark}

\section{Discussion}
\subsection{Summary}
In this work, we introduced a linear causal model for $P$-stationary processes, i.e., non-stationary discrete-time processes having all their $P$-periodic subsamplings jointly stationary. The associated graph contains a unique vertex for every $P$-periodic subsampling of each component process. A directed link between two subsamplings of two different component processes encodes an immediate influence, and this influence is quantified by the associated link filter, which describes how the past of the subsampled driving process linearly influences the subsampled target process. In addition, every process index has an associated $P$-stationary noise process, which we called the internal dynamics. We assumed the internal dynamics of two distinct components to be independent of each other. If the link filters are sufficiently small, the data of an $P$-stationary SEP give rise to a unique $P$-stationary observational process. 
Then we showed that structural vector autoregressive processes, whose coefficients are allowed to vary periodically over time with periodicity $P$, belong to the class of $P$-stationary structural equation processes by constructing an $P$-stationary SEP that is equivalent to the $P$-stationary SVAR process.

Our main interest was to infer the link filters in an $P$-stationary SEP associated with an $P$-stationary SVAR process when the observations are temporally incomplete, assuming that the summary graph of the SVAR process is known. The approach we took was based on three observations: First, every observational structure can be seen as a union of $RP$-periodic subsamplings, for appropriately selected $R$. Second, every $P$-stationary SEP has an equivalent representation as an $RP$-stationary SEP. Third, the problem can equivalently be solved in the frequency domain, where the algebraic relations are simpler. The first two observations allowed us to model an $P$-stationary SVAR processes with temporally incomplete observations as an $RP$-stationary SEP, where every missing subsample is treated as a latent process. However, we did not infer the link filters on the embedded $RP$-stationary SEP. Instead, we used the third observation to derive a linear relation between between the $P$-stationary spectrum and the $RP$-stationary spectrum, i.e., the embedded spectrum. Based on this relationship, we could determine which entries in the $P$-stationary spectrum can be recovered under the given observational structure. Hence, we arrived at the following reformulation of the initial problem: Is a link function identifiable from those entries in the $P$-stationary spectrum that are identifiable under the given observational structure? Based on this reformulation, we could derive a simple identification criterion for an instrumental process setting.

\subsection{Future work}
There are several lines in which the work presented in this paper can be extended. Instrumental processes provide a simple way to identify causal effects in latent process settings. However, in the literature on structural equation modeling, i.e., the case of i.i.d. observations, several other techniques have been explored to identify effects when the latent variables confound the observed ones in more interesting ways. In future work, one could try to identify link functions when the missing observations introduce more challenges. In this analysis it might be worthwhile to take into account that temporally incomplete observations hide entries in the spectrum in a slightly different way than when entire processes are unobserved. For example, to identify the spectral dependence between two processes, it suffices that only one of them is fully observable while observations of the other can be arbitrarily rare. However, when a process is not fully observed, the corresponding entry on the diagonal of the spectrum is not identifiable, no matter how rarely observations are missing.   

In the current paper, we only used the spectrum of the process, i.e., the cumulants of order two for the identification of causal effects. In ongoing work, we consider the case where the internal dynamics are filtered non-Gaussian white noise processes. Inspired by recent results in causal inference for linear non-Gaussian systems, we imagine that these ideas could also be applied to frequency domain versions of the higher order cumulants for the identification of the causal structure, i.e., the underlying directed graph, as well as the link functions, under more challenging observational structures. 

Finally, in this work we assumed to know the periodicity $P$ with which the coefficients in the underlying SVAR process are allowed to vary. In practice, this number $P$ would have to be estimated. A possibility for future work could be to derive an estimator and an associated test for $P$ based on the linear relationship between the embedded $RP$-stationary spectrum and the $P$-stationary spectrum.

\paragraph{Funding}
This project has received funding from the European Research Council (ERC) under the European Union’s Horizon 2020 research and innovation programme (grant agreement No 883818).

\paragraph{Statement on AI use}
The author consulted AI (Claude) when looking for related research, but for nothing else.

\appendix

\section{Regarding the embedding operator}
We will use the following two Lemmas to prove Lemma \ref{lemma: linear recovery of the embedded filter} and Lemma \ref{lemma: (tilted) convolution under the embedding} and Proposition \ref{prop: embedding properties}. 

\begin{lemma}
    Let $\lambda, \mu \in \ell^1(\Z; \C)$ and $i \in \Z$, then for convolution, it holds that 
    \begin{align*}
        \sh^i (\lambda \ast \mu ) = \sh^i(\lambda) \ast  \mu = \lambda \ast \sh^i (\mu),  
    \end{align*}
    and for tilted convolution, it holds that 
    \begin{align*}
        \sh^i (\lambda \ast' \mu) = \sh^i (\lambda ) \ast' \mu =\lambda \ast' \sh^{-i} (\mu).  
    \end{align*}
\end{lemma}
\begin{proof}
    We fix an integer $t\in \Z $, and note that $\sh^i(\lambda \ast \mu)(t) = \lambda \ast \mu (t-i)$. We now show the identities by direct computation:
    \begin{align*}
         \lambda\ast \mu (t-i )= \sum_{s\in \Z} \lambda (s) \mu(t-s-i ) = \sum_{s\in \Z } \lambda(s)\sh^i (\mu)(t-s) = \lambda \ast \sh^i (\mu ) (t). 
    \end{align*}
    Furthermore, using the index transformation $ s= s' -i$ on the right hand side of the first equation, we get that 
    \begin{align*}
        \lambda \ast \mu (t-i ) = \sum_{s' \in \Z } \lambda (s'-i) \mu (t-s') = \sh^i (\lambda) \ast \mu (t).  
    \end{align*}
    Regarding tilted convolution, we compute 
    \begin{align*}
        \lambda \ast' \mu(t-i) = \sum_{s \in \Z} \lambda(s+ t-i ) \mu (s) = \sh^i(\lambda) \ast' \mu(t),  
    \end{align*}
    and using the index shit $s = s' + i $ on the right hand side of the second equation, we get that 
    \begin{align*}
        \lambda \ast' \mu (t-i) = \sum_{s' \in \Z} \lambda(s'+t) \mu (s ' +i) = \lambda \ast' \sh^{-i}(\mu) ( t). 
    \end{align*}
\end{proof}

\begin{lemma}
    Let $\lambda, \mu \in \ell^1(\Z; \C)$ and let $R >0$ be a positive integer, then it holds that 
    \begin{align*}
        \defl_R (\lambda \ast \mu) = \sum_{m= 0}^{R-1} \defl_R \circ \sh^{-m}(\lambda)  \ast \defl_R \circ  \sh^m (\mu ), 
    \end{align*}
    and it also holds that
    \begin{align*}
        \defl_R (\lambda \ast' \mu ) &= \sum_{m = 0 }^{R-1} \defl_R \circ \sh^{-m} (\lambda) \ast' \defl_R \circ \sh^{-m } (\mu) . 
    \end{align*}
\end{lemma}
\begin{proof}
    For $t \in \Z$, we show the claim by direct computation: 
    \begin{align*}
        \lambda \ast \mu ( tR) &= \sum_{s\in \Z}  \lambda (s ) \mu (tR -s )  \\
        &= \sum_{m=0}^{R-1}\sum_{s' \in \Z } \lambda(s'R+m ) \mu ((t -s')R -m ) \\
        &= \sum_{m= 0}^{R-1} \sum_{s' \in \Z }  \defl_R \circ \sh^{-m} (\lambda )(s')  \defl_R \circ \sh^m(\mu)(t-s')  \\
        &= \sum_{m= 0}^{R-1}  \defl_R \circ \sh^{-m}(\lambda) \ast \defl_R \circ \sh^m (\mu ) (t).
    \end{align*}
    Now consider the tilted convolution
    \begin{align*}
        \lambda \ast' \mu (tR) &= \sum_{s \in \Z } \lambda( s +tR) \mu (s) \\
        &= \sum_{m = 0}^{R-1} \sum_{s' \in \Z } \lambda(m + (s'+t) R)\mu(s'R +m ) \\
        &= \sum_{m = 0}^{R-1}  \defl_R \circ \sh^{-m} (\lambda)\ast' \defl_R \circ \sh^{-m} (\mu) (t). 
    \end{align*}
\end{proof}

\begin{proof}[Proof of Lemma \ref{lemma: linear recovery of the embedded filter}]
    We need to show Equation (\ref{shifted recovery property}). Let us define $\tilde \lambda \coloneqq \mathsf{e}_R(\lambda)$, and let $\tilde \lambda_{i,j} \coloneqq \defl_R \circ \sh^{i-j}(\lambda)$. First, note that for $i,j \in [0,R-2]$, it holds that 
    \begin{align*}
        \tilde \lambda_{i,j} &= \tilde \lambda_{i+1, j+1}, & \tilde \lambda_{R-1, j-1} =  \sh^1 (\tilde\lambda_{0,j}), 
    \end{align*}
    and note that for any sequence $\lambda$ it holds that 
    \begin{align*}
        \left(\sh^{j} \circ \str_R(\lambda)\right) (k) &= \str_R(\lambda)(k-j) = \mathbf{1}_{R \Z}(k-j)\lambda\left(\mathbf{1}_{R\Z}(k-j)\frac{k-j}{R}\right).
    \end{align*}
    In particular, 
    \begin{align*}
        \left(\sh^{j} \circ \str_R \left( \tilde \lambda_{0,j} \right)\right)(k) &= \mathbf{1}_{R\Z}(k-j) \tilde \lambda_{0,j}\left( \frac{k-j}{R}\right) \\
        &= \mathbf{1}_{R\Z}(k-j) (\defl_R \circ \sh^{-j} (\lambda)) \left(\mathbf{1}_{R\Z}(k-j) \frac{k-j}{R} \right) \\
        &= \mathbf{1}_{R\Z}(k-j) \lambda (k).
    \end{align*}
    This implies the claimed identity for $i= 0$, i.e.,
    \begin{align}\label{decomposition: case i = 0}
        \mu = \sum_{j=0}^{R-1} \sh^j \circ \str_R \left( \tilde \lambda_{0,j} \right). 
    \end{align}
    Before we proceed with $i > 0$, we note that for any $i, j,R \in \Z $ with $R> 0$, the following relations hold:  
    \begin{align}\label{composing shift and stretch operators}
        \sh^{i} \circ \sh^{j} &= \sh^{i+j}, & \str_R \circ \sh^j &= \sh^{jR} \circ \str_R. 
    \end{align}
    Furthermore, note that $\sh^j$ and $\str_R$ are linear operators on the space of filters.  
    Using the two identities in (\ref{composing shift and stretch operators}), we proceed with $i > 0$ by direct computation
    \begin{align*}
        \sum_{j= 0}^{R-1} \sh^{j} \circ \str_R \left( \tilde\lambda_{i,j} \right) &= \sum_{j= 0}^{i-1} \sh^{j} \circ \str_R\left(\sh \left(\tilde\lambda_{0, R-i+j}\right)\right) + \sum_{j=i}^{R-1} \sh^{j} \circ \str_R\left( \tilde\lambda_{0,j-i} \right) \\ 
        &= \sum_{j= 0}^{i-1} \sh^{j+R} \circ \str_R \left( \tilde\lambda _{0, R- i+j}\right) + \sum_{j= i }^{R-1} \sh^{i+ (j-i)} \circ \str_R\left(\tilde\lambda_{0,j-i}\right) \\ 
        &= \sh^{i} \left(\sum_{j=0}^{i-1} \sh^{R-i+j} \circ \str_R \left( \tilde\lambda_{0, R-i+j} \right) + \sum_{j=i}^{R-1} \sh^{j-i} \circ\str_R \left( \tilde\lambda_{0,j-i} \right)\right) \\ 
        &= \sh^{i} \left( \tilde\lambda \right).
    \end{align*}
    For the last equality, we used (\ref{decomposition: case i = 0}) and the index transformation $j '= j + R -i$ for the first sum and $j ' = j +i $ for the second sum. 
\end{proof}
\begin{proof}[Proof of Lemma \ref{lemma: (tilted) convolution under the embedding}]
We begin with the identity for convolution. Let us define $\tilde \lambda \coloneqq \mathsf{e}_R(\lambda)$, and $\tilde \mu =\mathsf{e}_R (\mu )$ and $\tilde \gamma = \mathsf{e}_R (\gamma \coloneqq  \lambda \ast \mu )$. Furthermore, we set $\tilde \lambda_{i,j} = \defl_R \circ \sh^{i-j} (\lambda)$ and $\tilde \mu_{i,j} = \defl_R \circ \sh^{i-j}(\mu)$ and $\tilde \gamma_{i,j} = \defl_R \circ  \sh^{i-j} (\gamma)$ as the $i,j$-th block entry in $\tilde \lambda$ resp. $\tilde \mu$ and $\tilde \gamma$. We show the claimed identity by direct computation: 
\begin{align*}
    \tilde\gamma_{ij} &= \defl_R \circ \sh^{i-j} (\lambda \ast \mu) \\
    &= \defl_R ( \sh^{i}(\lambda) \ast \sh^{-j}(\mu) )\\
    &=\sum_{m = 0}^{R-1} \defl_R \circ \sh^{-m} (\sh^i (\lambda)) \ast \defl_R \circ \sh^m (\sh^{-j}(\mu)) \\
    &= \sum_{m = 0}^{R-1} \tilde \lambda_{i,m} \ast \tilde \mu_{m,j}, 
\end{align*}
which shows that $\tilde \gamma = \tilde \lambda \ast \tilde \mu$. 

Let us now proceed with the case of tilted convolution. For that purpose, we set $\tilde \mu' = \mathsf{e}_R' ( \mu)$ with $\tilde \mu'_{i,j}= \defl_R \circ \sh^{j-i} (\mu)$ and $\tilde \gamma' = \mathsf{e}_R (\gamma' \coloneqq \lambda \ast ' \mu )$ with $\tilde \gamma'_{i,j} = \defl_R \circ \sh^{i-j} (\gamma' )$. We show the claimed identity by direct computation: 
\begin{align*}
    \tilde \gamma_{i,j}' &= \defl_R \circ \sh^{i-j}(\lambda \ast'\mu) \\
    &= \defl_R( \sh^{i}(\lambda) \ast' \sh^j (\mu))  \\
    &=\sum_{m =0}^{R-1} \defl_R \circ \sh^{-m} (\sh^i(\lambda)) \ast' \defl_R\circ \sh^{-m} (\sh^j (\mu)) \\
    &= \sum_{m = 0}^{R-1} \tilde \lambda_{i,m} \ast' \tilde\mu'_{m,j}.
\end{align*}
The identity for embedding the filtered process follows by the same arguments as the identity for embedding the convolution of two sequences. To avoid repetition, we will not detail the verification.   
\end{proof}
\begin{proof}[Proof of Proposition \ref{prop: embedding properties}]
    The first three identities follow directly from Lemma \ref{lemma: (tilted) convolution under the embedding}. 
    So we only need to show that the computation of the ACS is compatible with the embedding. We denote the $P$-stationary ACS of $\X$ by $\Sigma$. Let $v',v \in V$ be two indices, then we compute for $t,k \in \Z$ as follows:  
    \begin{align*}
        [\mathrm{ACS}(\mathsf{E}_{R}(\X ))]_{v',v} (k) &= \mathbb{E}\left[ \mathsf{e}_{R}(\X_{v'}(t)) \mathsf{e}_{R}(\X_v(t-k ))^\top\right] \\
        &= \left[
            \E[\left(\X_{v'}(-i + tR) \right) \left(\X_{v}(-j+ (t-k)R)\right)^\top\right]_{i,j \in [0, R-1]} \\
            &= \left[ \Sigma_{v',v}(j-i + kR) \right]_{i,j \in [0, R-1]} \\
            &= \mathsf{e}_{R}(\Sigma_{v',v})(k). 
    \end{align*}
\end{proof}
\section{Representing a $P$-stationary SVAR representation as $P$-stationary SEP}
\subsection{Notation for time series graphs}
Let $\G = (V \times \Z, \D )$ be a $P$-stationary time series graph and let $\Phi \in \R^\G$ be a compatible filter. For $v,w \in V$ and $m,n \in [0,P-1]$, we write $\paths_{v,w}^{(m,n)}(k)$ to denote the set of directed paths on $\G$ starting at $v((t-k)P - m)$ and ending at $w(tP - n)$. For a directed path  
    \begin{align*}
        \pi &= (v_i(s_iP -l_i ) \to v_{i+1}(s_{i+1}P - l_{i+1}) \mid i \in [N-1]), 
    \end{align*}
    define the associated path filter as follows: 
    \begin{align*}
        \phi^{(\pi)} \coloneqq \left( \prod_{i=1}^{N-1} \phi^{(l_{i}, l_{i+1})}_{{v_i, v_{i+1}}}(s_{i+1}- s_{i}) \right) \mathbf{1}_{\{k\}}(j).
    \end{align*}
Given two directed paths, $\pi \in \paths_{u,v}^{(m,n)}(k)$ and $\pi' \in \paths_{v,w}^{(n,q)}(k')$, it holds that 
    \begin{align}\label{eq: convolution and concatenation tsg}
        \phi^{(\pi \ast \pi')} &= \phi^{(\pi)} \ast \phi^{(\pi')},  
    \end{align}
    where $\pi \ast \pi' \in \paths_{u,w}^{(m,q)}(k+k')$ is the concatenation of $\pi$ and $\pi'$. Finally, we define the union 
    \begin{align*}
        \paths_{v,w}^{(m,n)} &\coloneqq \bigcup_{k \in \Z} \paths^{(m,n)}_{v,w}(k), 
    \end{align*}
and associated with this set, we define 
\begin{align}\label{tsg path filter}
    T^{(m,n)}_{v,w}(\Phi) &\coloneqq \sum_{\pi \in \paths_{v,w}^{(m,n)}}  \phi^{(\pi)} = \sum_{k \in \Z } \sum_{\pi \in \paths^{(m,n)}_{v,w}(k)}\phi^{(\pi)}, 
\end{align}
which we arrange into the $P$-stationary sequence $T(\Phi): \Z \to \bmat(\R, V,[0,P-1])$. Note that $T(\Phi)$ is well-defined if the spectral radius of the matrix $\Phi(0)$ is less than one. 
\begin{lemma}\label{lemma: path sum tsg identity}
    It holds that 
    \begin{align*}
        \mathbb{I} + \Phi \ast T(\Phi ) &= T(\Phi ).
    \end{align*}
\end{lemma}
Before proving the lemma, we recall a basic property about directed paths. In order to state it, we define for $(v, m)\in V \times [0,P-1]$ the set 
\begin{align*}
    C(v,m) \coloneqq \left\{ (u,q) \in V \times [0,P-1] \mid \exists k:   v \to^{(m,q)}_k u\right\}.  
\end{align*}
\begin{lemma}\label{lemma: tsg path bijection}
For $m,n \in [0,P-1]$ and $v,w \in V$, there is the following bijection of sets 
    \begin{align*}
        \paths_{v,w}^{(m,n)} &\cong \begin{cases}
            \bigcup_{(u,q) \in C(v,m))} \paths_{u,w}^{(q,n)} & (v,m) \neq (w,n), \\
            \bigcup_{(u,q) \in C(v,m)} \paths_{u,w}^{(q,n)} \cup \{\epsilon_{(v,m)} \} & (v,m)=(w,n). 
        \end{cases}
    \end{align*}  
\end{lemma}
\begin{proof}
    Let $$\pi = (v_i (s_i P - l_i ) \to v_{i+1}(s_{i+1} P - l_{i+1})\mid i \in [N]) \in \paths_{v,w}^{(m,n)},$$
    then the claimed bijection maps $\pi$ to the path one obtains by removing the first link. Note that this map takes the empty path to the empty path. It is clear that this gives a bijection. 
\end{proof}
\begin{proof}[Proof of Lemma \ref{lemma: path sum tsg identity}]
    We prove the claimed identity element-wise and we abbreviate $T(\Phi)$ by $T$. Let therefore $v,w \in V$ and $m,n \in [0,P-1]$, then 
    \begin{align*}
        \mathbf{1}_{\{(v,m)\}}(w,n) + [\Phi \ast T]^{(m,n)}_{v,w} &= \mathbf{1}_{\{(v,m)\}}(w,n) +\sum_{u \in V} \sum_{q \in [0,P-1]} \Phi^{(m,q)}_{v,u} \ast T^{(q,n)}_{u,w} \\
        &= \mathbf{1}_{\{(v,m)\}}(w,n)+\sum_{(u,q) \in C(v,m)} \Phi^{(m,q)}_{v,u} \ast T^{(q,n)}_{u,w} \\
        &= T^{(m,n)}_{v,w},
    \end{align*}
    where the second equality follows from compatibility of $\Phi$ with $\G$ and the third equality follows from Lemma \ref{lemma: tsg path bijection} and Equation (\ref{eq: convolution and concatenation tsg}). 
\end{proof}
\begin{proposition}\label{proposition: pfd tsg}
    For $\Phi \in \R^\G$ such that the spectral radius of $\Phi(0)$ is less than one, it holds that 
    \begin{align}\label{Power-Psi: PFD}
        S(\Phi) &= T(\Phi). 
    \end{align}
\end{proposition}
\begin{lemma}\label{lemma: filter uniqueness}
    Let $\Lambda: \Z \to \R^{N \times N}$ be filter such that $I -\Lambda(0)$ is invertible and $\Lambda(k)= 0$ for $k<0$, and $A: \Z \to \R^{N \times N}$. Then the following implication holds:  
    \begin{align*}
       A&= \Lambda \ast A \implies A =0.
    \end{align*}
\end{lemma}
\begin{proof}
We show the claim by induction. First, consider $k=0$, then 
\begin{align*}
    (\Lambda \ast A) (0) = A(0) \iff (I - \Lambda(0))A(0) = 0 \iff A(0) =0,
\end{align*}
where for the first equality we used that $\Lambda(k)=0$ for $k< 0$, and for the last equivalence we used that $I - \Lambda(0)$ is invertible. 

Assuming $A(0)=A(\pm 1) \dots = A(\pm k) = 0$, we compute 
\begin{align*}
    A(k+1) &= \sum_{l=0}^{k+1} \Lambda(l)A(k+1-l) = \Lambda(0) A(k+1),
\end{align*}
and therefore  
\begin{align*}
     0 = (I -\Lambda(0)) A(k+1) \iff A(k+1) = 0.
\end{align*}
The case for $A(-k-1)$ works analogously. 
\end{proof}
\begin{proof}[Proof of Proposition \ref{proposition: pfd tsg}]
    The claim follows from the following identities: 
    \begin{align*}
        S(\Phi)& = \mathbb{I} + \Phi \ast S(\Phi ), & T(\Phi ) &= \mathbb{I } + \Phi \ast T(\Phi),
    \end{align*}
    which yield the equality $$ S(\Phi)- T(\Phi) = \Phi \ast (S(\Phi)- T(\Phi)),$$
    which implies, by Lemma \ref{lemma: filter uniqueness}, that $S(\Phi)=T(\Phi)$. 
\end{proof}
\subsection{The internal and link filters as sums of path filters}
Let $v,w \in V$ and $m,n \in [0,P-1]$ be such that $(v,m) \to (w,n)$ is a link on the process graph, then write $\paths_{v\to w}^{(m,n)}(k) \subset \paths_{v,w}^{(m,n)}(k )$ for the set of directed paths on $\G $ of the from
\begin{align*}
    \pi &= (v_i(s_iP -l_i) \to v_{i+1}(s_{i+1}P - l_{i+1}) \mid i \in [N-1]),
\end{align*}
where $v_i = w$ for $i > 1$ and $l_1 = m$ and $l_N = n$ and $s_N - s_1 = k$, and define 
\begin{align*}
    \paths^{(m,n)}_{v\to w} &\coloneqq \bigcup_{k \in \Z } \paths^{(m,n)}_{v\to w}(k ).
\end{align*}
The following is a special instance of Proposition \ref{proposition: pfd tsg}. 
\begin{proposition}\label{proposition: tsg decomposition sep filters}
For $v\neq w \in V$ and $m,n \in [0,P-1]$, the associated internal resp. link filters can be described as sums of TSG path filters in the following way:
    \begin{align*}
        \iota_v^{(m,n)} &= \sum_{\pi \in \paths_{v \to v}^{(m,n)}} \phi^{(\pi)}, & \lambda^{(m,n)}_{v,w} &= \sum_{\pi \in \paths^{(m,n)}_{v\to w}} \phi^{(\pi)}
    \end{align*}
    In particular, $v^\internal_m \not\to v_n \implies \iota^{(m,n)}_{v} = 0$, and $v_m \not\to w_n \implies \lambda^{(m,n)}_{v,w} = 0$. 
\end{proposition}
\begin{proof}
    Let $\G_v$ be the restriction of $\G$ to the vertex set $\{ v \} \times \Z$, and let $\Phi_v \in \R^{\G_v}$ be the restriction of $\Phi$ to $\G_v$, then identity regarding the internal filter follows by applying Proposition \ref{proposition: pfd tsg} to $\G_v$ and $\Phi_v$. The identity regarding the link filters follows then directly from the identity for the internal filters and the definition of the link filters, i.e., $\lambda_{v,w} = \Phi_{v,w} \ast \iota_w$. 
\end{proof}

\subsection{Path-wise decomposition of the power sum of the link filter}
Let $G = (V \times [0,P-1], D)$ be the process graph associated with the time series graph $\G$. We denote by $\paths((v,m), (w,n))$ the set of directed paths on $G$ starting at $(v,m)$ and ending at $(w,n)$. For a directed path $\pi = ((v_i, l_i) \to (v_{i+1}, l_{i+1})\mid i \in [N-1])$ on $G$, we define its associated path filter 
    \begin{align*}
        \lambda^{(\pi)} &\coloneqq \lambda_{v_1, v_2}^{(l_1, l_{2})} \ast \dots \ast \lambda_{v_{N-1}, v_N}^{(l_{N-1}, l_N )}.
    \end{align*}
For two directed paths, $\pi  \in \paths((u,m), (v,n))$ and $\pi' \in \paths((v,n), (w,q))$, it holds that 
    \begin{align*}
        \lambda^{(\pi \ast  \pi ' )} &= \lambda^{(\pi)} \ast \lambda^{(\pi')}, 
    \end{align*}
    where $\pi \ast \pi' \in \paths((u, m), (w,q))$ is the concatenation of $\pi$ and $\pi'$. 
\begin{lemma}\label{lemma: power lambda as sum of path filters}
    It holds that
    \begin{align*}
        S(\Lambda)_{v,w}^{(m,n)} &= \sum_{\pi \in \paths((v,m), (w,n))} \lambda^{(\pi)}. 
    \end{align*}
\end{lemma}
\begin{proof}
    For $k\geq 0$, define $\paths^k \subset  \paths((v,m),(w,n))$ to be the subset of directed paths $\pi$ consisting of exactly $k$ links. Using the fact that $\Lambda$ is compatible with $G$, one can show by induction that 
    \begin{align*}
        [\Lambda^k]^{(m,n)}_{v,w} &= \sum_{\pi \in \paths^k} \lambda^{(\pi)}.
    \end{align*}
    Then the claim follows as $$\paths((v,m), (w,n)) = \bigcup_{k \geq 0} \paths^k.$$
\end{proof}
To describe $S(\Lambda)$ in terms of path filters on $\G$, we introduce the following sets of directed paths. For $v,w \in V$ and $ m,n \in [0,P-1]$, define $\hat\paths_{v,w}^{(m,n)} \subset \paths^{(m,n)}_{v,w}$ to be the set of directed paths $\pi$ of the form $\pi = \pi' \ast \pi''$, where $\pi' \in \paths^{(m,n')}_{v \to v' }$ and $\pi'' \in \paths_{v',w}^{(n',n)}$, for some $v'\in V\setminus \{v\}$ and $n' \in [0,P-1]$.
\begin{lemma}
    It holds that 
    \begin{align}\label{eq: power filter process graph as sum of tsg path filters}
        S(\Lambda)^{(m,n)}_{v,w} &= \sum_{\pi \in \hat\paths^{(m,n)}_{v,w}} \phi^{(\pi)}. 
    \end{align}
\end{lemma}
\begin{proof}
    For a path $\pi = ((v_i, l_i) \to (v_{i+1}, l_{i+1}) | i \in [N-1]) \in \paths((v,m)), (w,n)$, we define 
    $$\paths^\pi \coloneqq \prod_{i\in [N-1]} \paths_{v_i \to v_{i+1}}^{(l_i, l_{i+1})}.$$ Note that concatenation of paths on $\G$ yields the injective map $$\ast:\paths^\pi \longrightarrow \hat\paths^{(m,n)}_{v,w}, \quad (\pi_i \mid i \in [N-1]) \mapsto \pi_1 \ast \cdots \ast \pi_{N-1}.$$
    Furthermore, Proposition \ref{proposition: tsg decomposition sep filters} in combination with the definition of $\lambda^{(\pi)}$ and $\paths^\pi$ gives that
    \begin{align}\label{eq proof: path filter process graph}
        \lambda^{(\pi)} = \sum_{(\pi_i \mid i \in [N-1]) \in \paths^\pi} \underbrace{\phi^{(\pi_1 )} \ast \cdots \ast \phi^{(\pi_{N-1})}}_{\phi^{(\pi_1 \ast \cdots \ast \pi_{N-1})}}.
    \end{align}

    On the other hand, for every path $\hat\pi \in \hat\paths_{v,w}^{(m,n)}$, there is a unique $\pi \in \paths((v,m), (w,n))$ such that $\hat\pi \in \paths^\pi$. In particular, concatenation on $\G$ yields a bijection 
    \begin{align}\label{eq proof: bijejction path sets}
         \bigcup_{\pi \in \paths((v,m), (w,n))} \paths^\pi \overset{\cong}{\longrightarrow}  \hat\paths^{(m,n)}_{v,w}.
    \end{align}
    Combining the identity (\ref{eq proof: path filter process graph}) with the bijection (\ref{eq proof: bijejction path sets}) and Lemma \ref{lemma: power lambda as sum of path filters} proves the identity (\ref{eq: power filter process graph as sum of tsg path filters}). 
\end{proof}

\begin{proposition}\label{prop: Power Phi is iota conv Power Lambda}
    It holds that 
    \begin{align*}
        S(\Phi) &= \iota  \ast S(\Lambda). 
    \end{align*}
\end{proposition}
\begin{proof}
    The proposition follows immediateyl by combining Equation (\ref{eq: convolution and concatenation tsg}) with the observation that concatenation on $\G$ yields a bijection $$ \paths_{v,w}^{(m,n)} \cong  \bigcup_{m' \in [0,P-1]}\paths_{v,v}^{(m,m')} \times \hat\paths_{v,w}^{(m', n)}.$$
\end{proof}
\begin{proof}[Proof of Proposition \ref{prop: P-stationary SEP representation of SVAR process}]
The solution of the SVAR process is $S(\Phi)^\top \ast \eta$, and the solution of the associated SEP is $S(\Lambda)^\top  \ast \iota^\top \ast \eta $. Hence, the solutions coincide by Proposition \ref{prop: Power Phi is iota conv Power Lambda}.
\end{proof}   

\section{Stability of $P$-stationary SVAR processes}
Let $(\G, \Phi, \eta)$ be a $P$-stationary SVAR model and $(G, \Lambda, \X^\internal)$ the associated SEP. We wish to develop a criterion for the parameter $\Phi$ that, if satisfied, implies stability of $\Lambda$. For that purpose, we will develop an alternative description of $S(\Phi)$ based on which we develop a criterion for $\Phi$ to check if $S(\Phi)$ is absolutely summable. Then we will use Proposition \ref{prop: Power Phi is iota conv Power Lambda} to conclude that $\Lambda$ is stable. 

For the alternative formulation of $S(\Phi)$, we construct a $P$-stationary $\Psi: \Z \to \bmat(\R , V, [0,P-1]) $, which is compatible with the summary graph $\smry(\G)$.
The entry of $\Psi$ corresponding to $v,w \in V,W$ and $m,n \in [0,P-1]$ is defined as follows: 
\begin{align*}
    \Psi_{v,w}^{(m,n)} &\coloneqq \begin{cases}
        S\left(\Phi_{v,v}^{(m,m)}\right) &  m = n \text{ and } v = w, \\
        \Phi_{v,w}^{(m,n)} \ast \Psi^{(n,n)}_{w,w} & m \neq n \text{ or } v \neq w  .
    \end{cases}
\end{align*}
Note that $\Psi^{(m,n)}_{v,w}$ is non-zero only if $(v,m) \to (w,n)$ on the auxiliary graph $\smry(\G)$. In order to describe the entries of $\Psi$ as sums of path filters on the time series graph $\G$, we define $\tilde\paths^{(m,n)}_{v\to w}(k) \subset \paths_{v\to w}^{(m,n)}(k)$ to be the directed paths on $\G$ of the following form:
\begin{align*}
     \pi &= (v_{i}(s_i P - l_i) \to v_{i+1}(s_{i+1}P- l_{i+1}) \mid i \in [N-1]), 
\end{align*}
where $v_1 = v$ and $l_1 = m$, and $v_i= w$ and $l_i = n$ for every $i > 1$. By the same principle as in the previous section, we define $$\tilde \paths_{v\to w}^{(m,n)} \coloneqq \bigcup_{k \in \Z} \tilde \paths_{v \to w}^{(m,n)}(k) .$$
\begin{lemma}\label{lemma: auxiliary path filter decomposition of power psi}
    For $v,w \in V$ and $m,n \in [0,P-1]$, it holds that  
    \begin{align*}
        \Psi^{(m,n)}_{v,w} &=\sum_{\pi \in \tilde\paths^{(m,n)}_{v\to w}} \phi^{(\pi)}.
    \end{align*}
\end{lemma}
\begin{proof}
    To prove the claim, we inspect the filter $S= S(\Phi^{(m,m)}_{v,v})$. In fact, it follows from Proposition \ref{proposition: pfd tsg} that
    \begin{align*}
        S &= \sum_{\pi \in \tilde\paths_{v\to v}^{(m,m)}} \phi^{(\pi)}. 
    \end{align*}
    Specifically, this follows by applying Proposition \ref{proposition: pfd tsg} to the $P$-stationary filter $$\Phi'^{(p,q)} = \begin{cases}
        \Phi^{(m,m)}_{v,v} & (p,q) = (m,m), \\
        0 & \text{otherwise},
    \end{cases}$$
    which is compatible with the time series graph $\G_v' = (\{v \} \times \Z, \D)$, where there is $v \to_{k}^{(p,q)} v$ if and only if $p =q = m$ and $v \to^{(m,m)}_k v$ is a lagged relation on $\G$, and by the very definition of $\Psi$, this observation suffices to show the statement of the lemma. 
    
\end{proof}
Let $\pi = ((v_i, l_i) \to (v_{i+1}, l_{i+1}) \mid i \in [N-1]) \in \paths_{\smry(\G)}((v,m), (w,n))$, then define the associated path filter
\begin{align}\label{eq: path filter psi tsg}
    \Psi^{(\pi)}&= \Psi^{(l_1,l_1)}_{v_1, v_1 } \ast\Psi^{(l_1, l_2)}_{v_1, v_2} \ast \dots \ast \Psi^{(l_{N-1}, l_N)}_{v_{N-1}, v_N} =\sum_{\tilde\pi \in \tilde\paths^\pi} \phi^{(\tilde\pi)}, 
\end{align}
where $\tilde\paths^{\pi} = \paths_{v_1, v_1}^{(l_1, l_1,)}\times \prod_{ i \in [N-1]} \tilde\paths_{v_i \to v_{i +1}}^{(l_i, l_{i+1})}$, and $$\phi^{(\tilde\pi)} = \phi^{(\tilde\pi_0)} \ast \dots \ast\phi^{(\tilde\pi_{N-1})}.$$

\begin{proposition}\label{prop: power psi and power phi equality}
    Given the notation as above, it holds that 
    \begin{align*}
        S(\Psi) &= S(\Phi). 
    \end{align*}
\end{proposition}
\begin{proof}
    Combining the fact that 
    \begin{align*}
        S(\Psi)_{v,w}^{(m,n)} &= \sum_{\pi \in \paths_{\smry(\G)}((v,m), (w,n))} \Psi^{(\pi)},
    \end{align*}
    with the observation that concatenation of paths on $\G$ induces a bijection of sets 
    \begin{align*}
        \bigcup_{\pi \in \paths((v,m), (w,n))} \tilde\paths^\pi &\cong \paths_{v,w}^{(m,n )}, 
    \end{align*}
    and the Equation (\ref{eq: path filter psi tsg}), we get the statement of the proposition
    \begin{align*}
        S(\Psi)^{(m,n)}_{v,w}=\sum_{\pi \in \paths_{\smry(\G)}((v,m), (w,n))}\psi^{(\pi)}&= \sum_{\pi \in \paths^{(m,n)}_{v,w}} \phi^{(\pi)} = S(\Phi)^{(m,n)}_{v,w}
    \end{align*}
\end{proof}
In the following, $\tilde V_1 \dots , \tilde V_N \subset V \times [0,P-1]$ denote the strongly connected components of the summary graph $\smry(\G)$, i.e., $(v,m), (w,n) \in \tilde V_i$ for some $i \in [N]$ if and only if $(v,m) \in \anc_{\smry(\G)} (w,n) \cap \des_{\smry(\G)}(w,n)$.   
\begin{proposition}\label{prop: stability P-stationary processes}
    Let $\Phi \in \R^\G$ be such that $\lVert [\Phi^\top]_{\tilde V_i, \tilde V_i } \rVert_1 < 1$ for all $i \in [N]$, then $S(\Phi)$ as well as $S(\Lambda)$ are absolutely summable. Furthermore, the Fourier transformed filter $\F(S(\Phi))$ takes the form of a rational function that is well-defined for every $z \in \C$ such that $|z|\leq 1$.
\end{proposition}
To alleviate the proof of Proposition \ref{prop: stability P-stationary processes}, we first prove a weaker version of it, which we state in the following lemma. For the proof of that lemma, it will be useful to recall the following fact about the operator norm with respect to the norm $\lVert \bullet \rVert_1$. Let $A \in \C^{V \times V}$ be a matrix, then it holds that 
\begin{align}\label{eq: operator one-norm}
    \lVert A \rVert_1 &= \max_{v \in \R^V \setminus \{0\}} \frac{\lVert A \rVert_1}{\lVert v \rVert_1} = \max_{i \in V} \sum_{j \in V} \lvert A_{i,j}\rvert. 
\end{align}
Note that $\lVert A^\top \rVert_1 < 1$ implies that the spectral radius of $A$ is less than one. Suppose $A v = \lambda v$, then $\lVert Av \rVert_1 = \lvert \lambda \rvert \lVert v \rVert_1$. Hence, $|\lambda| <1$ follows from 
\begin{align*}
    \lVert A v \rVert_1 = \sum_{i \in V} \lvert \sum_{j\in V} A_{i,j} v_j \rvert \leq \sum_{j \in V } \lvert v_j \rvert  \sum_{i \in V } \lvert A_{i,j}\rvert \leq \lVert v \rVert_1 \lVert A^\top \rVert_1 < \lVert v \rVert_1. 
\end{align*}
Recall from the main part that for a matrix-valued filter $\Phi$, we defined 
\begin{align*}
    \lVert \Phi \rVert_1 &\coloneqq \lVert \sum_{k \in \Z } \lvert \Phi(k) \rvert  \rVert_1 
\end{align*}
Furthermore, note that $\sum_{k \in \Z} \lvert \Phi(k)\rvert = \F(\lvert \Phi \rvert)(1)$, where $|\Phi|(k) = \lvert \Phi (k) \rvert $, for every $k\in \Z$.  
\begin{lemma}\label{lemma: stability P-stationary SVAR processes}
    If $\lVert \Phi^\top \rVert_1 < 1$, then the filter $S(\Phi)$ is absolutely summable, and $S(\F(\Phi))(z)$ takes the form of a rational function in $z$ that is well-defined for every $z \in \C$ with $|z|\leq 1$. 
\end{lemma}
\begin{proof}
    The preceding observation in combination with the assumption of the lemma give, for $|z|\leq 1$, the chain of inequalities $$\lVert\F(\Phi^\top)(z)\rVert_1 \leq\lVert\F(|\Phi|^\top)(1)\rVert_1 \leq \lVert \Phi^\top \rVert_1 < 1.$$
    In particular, the spectral radius of $\F(\Phi)(z)$ is less than one, for every $|z|\leq 1$, meaning that $ S(\F(\Phi)(z)) = S(\F(\Phi)(z)) = ( I - \F (\Phi)(z))^{-1}$ exists and takes the form of a rational function in $z$ because $\Phi$ being compatible with the finite order time series graph $\G$ implies that $\F(\Phi)$ is a polynomial in $z$. Finally, noting that $S(|\Phi|) = |S(\lvert \Phi \rvert )|$, we use the observations preceding the lemma to conclude with the absolute summability of $S(\Phi)$, which is implied by
    \begin{align*}
         \lVert S (\Phi) \rVert_1 \leq\lVert S(|\Phi|) \rVert_1 = \lVert \F(S(|\Phi|))(1) \rVert_1 = \lVert \sum_{k \geq 0 } \F(|\Phi|)^k(1) \rVert_1 < \infty.
    \end{align*}
\end{proof}
\begin{proof}[Proof of Proposition \ref{prop: stability P-stationary processes}]
    Since the assumption on $\Phi$ is stated with respect to the operator norm, which satisfies Equation (\ref{eq: operator one-norm}), we may assume without loss of generality that $\Phi$ is entry-wise non-negative, i.e., for all $k\in \Z$, every entry of the matrix $\Phi(k)$ is non-negative.
    
    Let $\tilde V_1, \dots , \tilde V_N$ be the strongly connected components of the summary graph $\smry(\G)$, and assume that they are ordered such that for every $(v,m)\in \tilde V_i$ and $(w,n)\in \tilde V_j $ with $i < j$, it holds that $\paths_{\smry(\G)} ((v,m), (w,n))= \emptyset$. 
    
    Now let $i \leq j$, and $(v,m)\in \tilde V_i$ and $(w,n)\in \tilde V_j$. Furthermore, let us denote by $I(i,j ) = \{ x_1 < \cdots < x_k \mid x_l \in [i,j], x_1 = i, x_k = j\}.$ Then note that concatenation of paths on $\smry(\G)$ yields an isomorphism of sets: 
    \begin{align}\label{eq: path set decomposition auxiliary graph}
        \paths_{\smry(\G)}(\tilde V_i , \tilde V_j)  \cong \bigcup_{ \mathbf{x}\in I(i,j ) } \paths_{\smry(\G)}(\tilde V_{x_1}, \tilde V_{x_1}) \times \prod_{l=2}^{k} D'_{\tilde V_{x_{l-1}}, \tilde V_{x_l} }\times \paths_{\smry(\G)}(\tilde V_{x_l}, \tilde V_{x_l }),   
    \end{align}
    where $D'_{\tilde V_x, \tilde V_y} = D \cap (\tilde V_x \times \tilde V_y)$, and $\paths_{\smry(\G)} (\tilde V_x, \tilde V_y ) = \bigcup_{(v,w )\in \tilde V_x \times \tilde V_y}^k \paths_{\smry(\G)}(v,w)$.

    By combining Lemma \ref{lemma: auxiliary path filter decomposition of power psi} with the bijection (\ref{eq: path set decomposition auxiliary graph}), we get that 
    \begin{align}\label{equation: decomposition power psi}
        [S(\Psi)]_{\tilde V_i , \tilde V_j} &= \sum_{i = x_1 <\dots < x_k=j} [S(\Psi)]_{\tilde V_{x_1}, \tilde V_{x_1}} \ast \prod_{l=2}^k [\Psi]_{\tilde V_{x_{l-1}}, \tilde V_{x_l}} \ast [S(\Psi)]_{\tilde V_{x_l}, \tilde V_{x_l}}.  
    \end{align}
    The right hand side is a finite sum in which every summand is  convolution of finitely many factors. In particular, to prove the proposition, it suffices to show that each of those factors is absolutely summable. In order to do that, note that the definition of strictly connected components and Proposition \ref{prop: power psi and power phi equality} imply the chain of identities 
    \begin{align*}
        [S(\Psi)]_{\tilde V_x, \tilde V_x } = S([\Psi]_{\tilde V_x, \tilde V_x }) = S([\Phi]_{\tilde V_x, \tilde V_x }),
    \end{align*}
    where the assumption on $\Phi$ together with Lemma \ref{lemma: stability P-stationary SVAR processes} yields that the filter on the right hand side is absolutely summable. To see that $\Psi$ is absolutely summable, note that the assumption on $\Phi$ implies that $\lVert \Phi_{v,v}^\top \rVert_1 < 1$. Hence, the definition of $\Psi$ in combination with Lemma \ref{lemma: stability P-stationary SVAR processes} (applied to an appropriate subgraph of $\G$) implies that $\Psi$ is absolutely summable. To summarize, the filter $S(\Psi)$ is absolutely summable. Finally, combining Proposition \ref{prop: power psi and power phi equality} with Proposition \ref{prop: Power Phi is iota conv Power Lambda} gives that $\iota \ast S(\Lambda) = S(\Phi ) = S(\Psi)$. Since we assumed $\Phi$ to be entry-wise non-negative, it follows that also $\Lambda$ and $\iota$ are entry-wise non-negative. Furthermore, since $\iota(0)$ is the identity matrix, it holds entry-wise that $0\leq S(\Lambda)(k) \leq S(\Phi)(k)$, for every $k \in \Z $. In particular, the absolute summability of $S(\Phi)$ implies that $S(\Lambda)$ is absolutely summable as well.  

    Finally, that $\F(S(\Phi))$ takes the form of a rational function that is well defined for every $|z| \leq1$ follows by applying the Fourier transform to Equality (\ref{equation: decomposition power psi}), as it results in a sum of products of rational functions, and all of those are well-defined for $|z| \leq 1$ by Lemma \ref{lemma: stability P-stationary SVAR processes}. 
\end{proof}
We conclude this section with the observation that $S(H)$ admits an explicit parameterization in terms of the SVAR filter $\Phi$. Associated with $\Phi$, we define $\Phi^\internal: \Z \to  \bmat(\R, V, [0,P-1])$ to be the $P$-stationary filter whose block entry $\Phi^\internal_{u,v} $ is equal to $ \Phi_{v,v}$ if $u =v$, and equal to the zero-filter otherwise.
\begin{corollary}
    Suppose $(\G, \Phi, \eta)$ is a $P$-stationary SVAR model, where $\Phi$ is as in Proposition \ref{prop: stability P-stationary processes}, then 
    \begin{align*}
         S(H) & = (I - \F(\Phi^\internal))(I - \F(\Phi ))^{-1}.
    \end{align*}
\end{corollary}

\section{The $P$-stationary SEP representation of a stationary SVAR process}
\begin{proof}[Proof of Proposition \ref{prop: field of special Toeplitz matrices}]
    First, it follows from the definition that $\mathcal{T}_r(R,P)$ is closed under addition and scalar multiplication. Let $\mathbf{A}, \mathbf{B} \in  \mathcal{T}_r(R,P)$.

    \textsc{Multiplicative:} We now show that $\mathbf{A}\mathbf{B} = \in \mathcal{T}_r(R,P)$. Let $A_0, \dots , A_{R-1} \in \R(z)^{P \times P }$ and $B_0, \dots , B_{R-1} \in \R(z)^{P\times P }$ the matrices specifying $\mathbf{A}$ resp. $\mathbf{B}$. Let us compute, for $i,j \in [0,P-1]$: 
    \begin{align*}
        C_j & \coloneqq [\mathbf{A}\mathbf{B}]_{0,j} = \sum_{l=0}^{R-1} A_{0,l} B_{l,j} = \sum_{l= 0 }^{j} A_l B_{j-l} + \sum_{l = j+1}^{R-1} r A_lB_{R-l+j}
    \end{align*}
    Let $i,j \in [0,R-1]$ and assume that $i \leq j$, we wish to show that $[AB]_{i,j} = \gamma_{j -i}$.  
    \begin{align*}
        [\mathbf{AB}]_{i,j} &= \sum_{l=0}^{i-1} A_{i,l}B_{l,j} + \sum_{l=i}^j A_{i,l}B_{l,j} + \sum_{l=j+1}^{R-1}A_{i,l}B_{l,j} \\
        &= \sum_{l=0}^{i-1} rA_{R-i +l} B_{j-l} + \sum_{l=i}^{j} A_{l-i} B_{j-l} + \sum_{l=j+1}^{R-1} A_{l-i}rB_{R-l+j}
    \end{align*}
    Now make the following substitutions: $l = - R + i +l'$ for the first summand, $l = l' +i $ for the other two summands, which yields: 
    \begin{align*}
        [\mathbf{AB}]_{i,j} &= \sum_{l'= R-i}^{ R-1} r A_{l'}B_{R-l' + (j-i)} + \sum_{l'=0}^{j-i}A_{l'} B_{(j-i)-l'} + \sum_{l'=j-i+1}^{R-i-1}A_{l'}r B_{R-l'+ (j-i)} \\
        &= \sum_{l'=0}^{j -i} A_{l'} B_{(j-i)-l'} + \sum_{l'= j -i +1}^{R-1} r A_{l'}B_{R- l' + (j-i)} \\
        &= C_{j-i}.
    \end{align*}
    Now assume $i > j$ and compute 
    \begin{align*}
        [\mathbf{AB}]_{i,j} &= \sum_{l=0}^{j} A_{i,l}B_{l,j} + \sum_{l=j+1}^{i-1} A_{i,l}B_{l,j} + \sum_{l=i}^{P-1}A_{i,l}B_{l,j} \\
        &= \sum_{l= 0}^j r A_{R-i+l}B_{j-l} + \sum_{l = j+1}^{i-1}r A_{R -i +l }rB_{R-l+j} + \sum_{l=i}^{R-1}A_{l-i} rB_{R-l+j} \\
        &= \sum_{l'= R-i}^{R-i+j} r A_{l'}B_{R-i+j-l'} + \sum_{l'= R -i+j+1}^{R-1}r A_{l'}r B_{2R-i+j-l'} + \sum_{l'= 0 }^{R-i-1}A_{l'}rB_{R-i+j -l'} \\
        &= r\left(\sum_{l'= 0}^{R - i+j }A_{l'} B_{R-i+j -l'} + \sum_{l'= R-i+j+1}^{R-1}r A_{l'} B_{2R-i+j-l'} \right) \\
        &= r C_{R-i+j}. 
    \end{align*}
    This shows that the set $\mathcal{T}^P_r $ is multiplicative. 
    

    \textsc{Closed under inversion:} Let $\mathbf{A} \in \mathcal{T}_r(R,P)$ be invertible over $\R(z)$, i.e., $\det(\mathbf{A}) \neq 0$ as a rational function. Let $p_\mathbf{A}(\lambda) = \det(\lambda I - \mathbf{A})$ be the characteristic polynmial over $\R(z)$ taking the form 
    \begin{align*}
        p_\mathbf{A}(\lambda) &= \lambda^N  + f_{N-1} \lambda^{N-1} + \cdots + f_1 \lambda + f_0, 
    \end{align*}
    where $N= RP$ and $f_i \in \R(z)$. Note that $f_0 = (-1)^{N}\det(\mathbf{A})$. The Cayley Hamilton Theorem now implies that $p_\mathbf{A}(\mathbf{A}) = 0$, which gives the following expression for the inverse
    \begin{align*}
        \mathbf{A}^{-1} &= (-1)^{N+1} \det(\mathbf{A})^{-1} (A^{N-1} + f_{N-1} A^{N-2}+ \cdots + f_2 A+ f_1 I). 
    \end{align*}
    Since $\mathcal{T}_r(R,P)$ is closed under multiplication, addition and scalar multiplication, it follows that $\mathbf{A}^{-1}$ is an element in $\mathcal{T}_r(R,P)$. 
\end{proof}

\section{Graphical criterion for non-vanishing minors}
\begin{proof}[Proof of Proposition \ref{prop: non-intersecting paths}]
    Let $T = (\tau_1, \dots, \tau_N) \in \treks_{\smry(\G)}(X, Y)$ be a non-intersecting trek system, where $X = ((v_i, m_i) \mid i \in [N]) $ and $Y = ((y_i, n_i \mid i \in [N]))$. Now define a time series subgraph $\tilde\G \leq \G $. First, we will show that $T$ can always be adapted such that every vertex visited by some unique $\tau_i$ is visited by $\tau_i$ exactly once. Note that $T$ being non-intersecting implies that every vertex is visited by at most one $\tau_i$. Let $v$ be a vertex appearing on trek $\tau_i = (\pi_i^R  , \trektop_i , \pi_i^L )$, where 
    \begin{align*}
        \pi_i^R &= \trektop_i \to v_1^R \to \dots \to v_{M_i}^R \\
        \pi_i^L &= \trektop_i \to v_1^L \to \dots \to v_{Q_i}^L.
    \end{align*}
    Then define $a_R =  \min \{ i\mid v_i^R = v \}$ and $a_L$ analogously. Furthermore, define $b_R = \max\{ i \mid v_i^R = v\} $ and $b_L$ analogously. If $v$ appears on only one side, say the left hand side, then replace $\tau_i$ by $\tau_i' = (\pi_i^{L'}, \trektop_i, \pi^R_i)$, where  
    \begin{align*}
        \pi^{L'}_i &= \begin{cases}
            \trektop_i \to v_1^L \to \dots \to  v^L_{a_L} \to v_{b_L+1}^L \to \dots \to v_{Q_i}^L & b_L < Q_i, \\ 
            \trektop_i \to v_1^L \to \dots \to  v^L_{a_L} & b_L = Q_i.
        \end{cases}  
    \end{align*}
    If $v$ appears on both sides, then we replace $\tau_i$ by $\tau_i' = (\pi_i^{L'}, v, \pi^{R'}_i)$, where 
    \begin{align*}
        \pi_i^{L'} &= v \to v_{b_L +1}^L \to \dots \to v_{Q_i}^L , \\
        \pi_i^{R'} &= v \to v_{b_R +1}^R \to \dots \to v_{M_i}^R 
    \end{align*}
    Note that the trek $\tau_i'$ visits the vertex $v$ exactly once. Proceeding this way, any non-intersecting trek system can be reduced to a trek system which visits every vertex at most once. So from now on we assume $T$ to be such a trek system. 
    
    For every link $(u,l) \to (u',l')$ appearing on some trek $\tau_i$, we add the relation $u \to_k ^{(l,l')} u'$ to $\smry(\G)$, where 
    \begin{align*}
        k = \min \{j \mid u\to_{j}^{(l',l)} u' \text{ on } \G \}. 
    \end{align*}
    Furthermore, let $\smry(\tilde\G)$ be the summary graph of $\tilde\G$. By construction and the assumption that $T$ is non-intersecting, it follows that $ \smry(\tilde\G)$ is an acyclic subgraph of $\smry(\G)$. Secondly, we have that
    \begin{align*}
        \treks_{\tilde \smry(\G)}((v,m), (w,n)) &= \begin{cases}
            \{ \tau_i  \} & (v,m) = (v_i, m_i) \text{ and } (w,n) = (w_{\alpha(i)}, n_{\alpha(i)}), \\ 
            \{ \epsilon_{(v,m)} \} & (v,m) = (w,n), \\
            \emptyset & \text{otherwise}. 
        \end{cases}
    \end{align*}
    In particular, any two vertices are connected by at most one trek. Applying Proposition \ref{prop: trek-rule} to the subgraph $\tilde{G}$, we get the following expression for the spectrum 
    \begin{align*}
        [\fcum_{\tilde\G}]_{X,Y} &= \left( \mathbf{1}_{\{\alpha(i)\}}(j) \fcum_{\tilde\G}^{(\tau_i)}) \mid i,j \in [N] \right) \in \K_{\tilde \G}^{X\times Y}. 
    \end{align*}
    By construction, we have that 
    \begin{align*}
        \fcum_{\tilde\G}^{(\tau_i )} &= \omega_{z_i}^{(k_i)}\F(\psi^{(\pi^L_i)}) (\F(\psi^{(\pi^R_i)}) )^\ast,
    \end{align*}
    where $\trektop_i = (z_i, k_i)$ and $\F\left(\Psi^{(\pi_i^L)}\right) = \prod_{i=1}^{K_i-1} \F\left(\Phi_{v^L_i, v_{i+1}^L}^{(l_i, l_{i+1})}\right)$ and analogously for $\F\left(\Psi^{(\pi^R_i)}\right)$. Note that all of those are non-zero polynomials in the indeterminates of the SVAR coefficients, so that 
    \begin{align*}
        \det \left([\fcum_{\tilde\G}]_{X,Y}\right) =(-1)^\Pi \prod_{i= 1}^N \fcum_{\tilde \G }^{(\tau_i)}
    \end{align*}
    is a product of non-zero rational functions, and therefore non-zero overall, and since this implies that also $\det([\fcum_\G]_{X,Y})\in \K_\G$ is non-zero, we conclude that $[\fcum_\G]_{X,Y} \in \K_\G^{X,Y}$ is invertible. 
\end{proof}

\bibliographystyle{plain}
\bibliography{references}

\end{document}